\documentclass[12pt]{amsart}
\usepackage{amsmath,amssymb,amsthm}

\usepackage{upgreek}

\usepackage{enumitem}

\usepackage{ifthen,verbatim}
\usepackage{mathrsfs}
\usepackage{color}
\usepackage[urlcolor=blue,colorlinks=true]{hyperref}
\usepackage{url}
\usepackage{pdfpages}
\usepackage[english]{babel}
\usepackage{verbatim,here}
\usepackage[T1]{fontenc}
\usepackage{floatflt,graphicx,graphics}
\usepackage{a4wide}

\graphicspath{{figures}}

\numberwithin{equation}{section}
\nonstopmode
\theoremstyle{plain}

\newtheorem{corollary}[equation]{Corollary}

\newtheorem{theorem}[equation]{Theorem}

\newtheorem{lemma}[equation]{Lemma}

\newtheorem{proposition}[equation]{Proposition}
\newtheorem{problem}[equation]{Problem}

\theoremstyle{definition}
\newtheorem{rem}[equation]{Remark}

\newtheorem{example}[equation]{Example}
\newtheorem{definition}[equation]{Definition}

{\qed\bigskip}

\newcounter{alphabet}

\newcommand{\be}{\begin{eqnarray}}
\newcommand{\ee}{\end{eqnarray}}
\newcommand{\ba}{\begin{array}}
\newcommand{\ea}{\end{array}}
\newcommand{\ben}{\begin{eqnarray*}}
\newcommand{\een}{\end{eqnarray*}}

\newcommand{{\tth}}{\mathrm{th}}
\newcommand{{\sh}}{\mathrm{sh}}
\newcommand{{\ch}}{\mathrm{ch}}

\newcommand{\sn}{\mathrm{sn}}

\renewcommand{\Im}{{\,\operatorname{Im}\,}}
\renewcommand{\Re}{{\,\operatorname{Re}\,}}

\newcommand {\M} {\mathsf{M}}

\newcounter{minutes}
\divide\time by 60
\newcounter{hours}
\multiply\time by 60
\addtocounter{minutes}{-\time}

\begin{document}

\bibliographystyle{amsplain}

\def\thefootnote{}
\footnotetext{
\texttt{\tiny File:~\jobname .tex,
          printed: \number\year-\number\month-\number\day,
          \thehours.\ifnum\theminutes<10{0}\fi\theminutes}
}
\makeatletter\def\thefootnote{\@arabic\c@footnote}\makeatother

\title[Conformal capacity of hedgehogs]{Conformal capacity of hedgehogs}
\author[Dimitrios Betsakos, Alexander Solynin and Matti Vuorinen]{Dimitrios Betsakos, Alexander Solynin and Matti Vuorinen}

\address{Department of Mathematics, Aristotle University of
Thessaloniki, GR-54124 Thessaloniki, Greece}

\email{betsakos@math.auth.gr}

\address{Department of Mathematics and Statistics, Texas Tech
University, Box 41042, Lubbock, Texas 79409, USA}
\email{alex.solynin@ttu.edu}

\address{Department of Mathematics and Statistics, FI-20014
University of Turku, Finland} %
\email{vuorinen@utu.fi}


\begin{abstract}
We discuss problems concerning the conformal condenser capacity of
``hedgehogs'', which are compact sets $E$ in the unit disk
$\mathbb{D}=\{z:\,|z|<1\}$ consisting of a central body $E_0$ that
is typically a smaller disk $\overline{\mathbb{D}}_r=\{z:\,|z|\le
r\}$, $0<r<1$, and several spikes $E_k$ that are compact sets
lying on radial intervals $I(\alpha_k)=\{te^{i\alpha_k}:\,0\le
t<1\}$. The main questions we are concerned with are the
following: (1) How does the  conformal capacity ${\rm cap}(E)$ of
$E=\cup_{k=0}^n E_k$ behave when the spikes $E_k$, $k=1,\ldots,n$,
move along the intervals $I(\alpha_k)$ toward the central body if
their hyperbolic lengths are preserved during the motion? (2) How
does the
 capacity ${\rm cap}(E)$ depend on the distribution of
angles between the spikes $E_k$? We prove several results related
to these questions and discuss methods of applying symmetrization
type transformations  to study the  capacity of hedgehogs. Several
open problems, including problems on the
 capacity of hedgehogs in the three-dimensional
hyperbolic space,  will also be suggested.
\end{abstract}

\maketitle

{\hspace{10.3cm}{{\Small{In memoriam Jukka Sarvas (1944-2021).}}}}


 \hypersetup{linkcolor=black}
\tableofcontents
\section*{Notation} %
\begin{itemize}
\item[$\bullet$] %
$\mathbb{C}$ - complex plane.
\item %
$\mathbb{D}=\{z\in \mathbb{C}:\, |z|<1\}$ - open unit disk
centered at $z=0$.
\item %
$\lambda_{\mathbb{D}}(z)$ - density of the hyperbolic metric in
$\mathbb{D}$. %
\item %
$\ell_{\mathbb{D}}(\gamma)$ - hyperbolic length of $\gamma$. %
\item %
$A_{\mathbb{D}}(E)$ - hyperbolic area of $E$. %
\item %
$p_{\mathbb{D}}(z_1,z_2)=\left|(z_1-z_2)/(1-z_1\overline{z}_2)\right|$
- pseudo-hyperbolic metric in $\mathbb{D}$.
\item %
$d_{\mathbb{D}}(z_1,z_2)=\log
\frac{1+p_{\mathbb{D}}(z_1,z_2)}{1-p_{\mathbb{D}}(z_1,z_2)}$ -
hyperbolic metric in $\mathbb{D}$. %
\item %
$[z_2,z_2]$ and $(z_1,z_2)$ - closed and open Euclidean intervals
with end points $z_1$ and $z_2$. %
\item %
$[z_1,z_2]_h$ and $(z_1,z_2)_h$ - closed and open hyperbolic
intervals with end points $z_1$ and $z_2$.
\item %
${\rm cap}(E)={\rm cap}(\mathbb{D},E)$ - conformal capacity of a
compact set $E\subset \mathbb{D}$. %
\item %
 ${\rm log.cap}(E)$ - logarithmic capacity of $E\subset
\mathbb{C}$. %
\item %
$\M(\Gamma)$ - modulus of the family of curves $\Gamma$.
\item %
 $\mathcal{R}_\gamma(E)$ - reflection of a set $E$ with
    respect to the hyperbolic geodesic $\gamma$.%
    \item %
 $\mathcal{P}_\gamma(E)$ - polarization of a set $E$ with
    respect to the hyperbolic geodesic $\gamma$. %
\item %
 $\mathcal{K}(k)$ and
$\mathcal{K}'(k)=\mathcal{K}(\sqrt{1-k^2})$ - complete elliptic
integrals of the first kind.
\end{itemize}

\section{Introduction} %
The main theme discussed in this paper is the dependence of the
condenser capacity ${\rm cap}(\mathbb{D}, E)$ on the geometric
structure and characteristics of a compact set $ E \subset
\mathbb{D},$ where $ \mathbb{D}$ is the unit disk
$\mathbb{D}=\{z:\,|z|<1\}$ in the complex plane $\mathbb{C}$.
 For brevity, we call  ${\rm cap}(\mathbb{D},E) $
the conformal capacity of $E$ or the capacity of $E$ and denote it
by  ${\rm cap}(E)$.

Due to the conformal invariance of the capacity it is natural to
equip  $\mathbb{D}$ with the hyperbolic metric. Indeed, very
recently it was shown in \cite{NV2} that isoperimetric
inequalities in hyperbolic metric yield simple upper and lower
bounds for the capacity  in the case when $E$ is a finite union of
hyperbolic disks. In the subsequent work \cite{NV1, NRV1, NRV2}
these ideas were developed further, and it was also pointed out
that similar ideas were also applied by F.W. Gehring \cite{G} and
R. K\"uhnau \cite{Ku} fifty years earlier.

 In most cases, we deal with compact sets $E=\cup_{k=0}^nE_k$ consisting
 of a \emph{central body } $E_0$, which can be absent, and \emph{spikes} $E_k$,
 $k=1,\ldots,n$, that are closed intervals or any compact sets lying on $n$ radial intervals
 $I(\alpha_k)$, where $I(\alpha)=\{te^{i\alpha}:\,0\le t<1\}$, $I=I(0)$.
 This type of compact sets have appeared in several research
 papers, for instance, in a recent paper \cite{IRZ} by J.-W. M. Van Ittersum, B. Ringeling, and W.
 Zudilin, where the term ``hedgehog'' was suggested for this shape of compact
 sets. Interestingly enough, estimates of the capacity and other
 characteristics of hedgehogs appeared to be useful in studies on
 the Mahler measure and Lehmer's problem. Beside the above
 mentioned work of three authors, hedgehog structures appeared in
 the paper \cite{Pritsker} by I.~Pritsker and, earlier, the same hedgehog structure appeared in
 \cite{Solynin1998}.

 The hyperbolic metric in $\mathbb{D}$ is defined by the element of length 
$$ 
\lambda_{\mathbb{D}}(z)\,|dz|=\frac{2 |dz|}{1 - |z|^2}. %
$$  
Then the hyperbolic length $\ell_{\mathbb{D}}(E)$ of a compact
subset $E$ of
a rectifiable curve can be calculated as %
$$ 
\ell_{\mathbb{D}}(E)=\int_E\frac{2 |dz|}{1 - |z|^2}. %
$$ 
Furthermore, the hyperbolic geodesics are circular arcs in
$\mathbb{D}$ that are orthogonal to the unit circle
$\mathbb{T}=\partial \mathbb{D}$ at their end points. The
hyperbolic distance between points $z_1$ and $z_2$ in
$\mathbb{D}$, that is equal to the hyperbolic length
$\ell_{\mathbb{D}}([z_1,z_2]_h)$ of the closed hyperbolic
interval $[z_1,z_2]_h$ joining these points, is given by %

\begin{equation} \label{Hyperbolic distance}
d_{\mathbb{D}}(z_1,z_2)=\log
\frac{1+p_{\mathbb{D}}(z_1,z_2)}{1-p_{\mathbb{D}}(z_1,z_2)}, %
\end{equation}
where $p_{\mathbb{D}}(z_1,z_2)$ is the pseudo-hyperbolic metric
defined as %
\begin{equation} \label{Pseudo-hyperbolic distance}
p_{\mathbb{D}}(z_1,z_2)=
\left|\frac{z_1-z_2}{1-z_1\overline{z}_2}\right|.
\end{equation} %
Everywhere below, $[z_1,z_2]_h$ and $(z_1,z_2)_h$ stand,
respectively, for the closed and open hyperbolic intervals with
end points $z_1$ and $z_2$. Similarly, notations $[z_2,z_2]$ and
$(z_1,z_2)$ will be reserved for the closed and open Euclidean
intervals with end points $z_1$ and $z_2$. If $z_1$ and $z_2$ lie
on the same diameter of $\mathbb{D}$  then, of course,
$[z_1,z_2]_h=[z_1,z_2]$ and $(z_1,z_2)_h=(z_1,z_2)$.

When $z_1=0$ and $z_2=r$, $0<r<1$, the hyperbolic length
$\uptau=\uptau(r)$ of the interval $[0,r]$ and its Euclidean
length $r=r(\uptau)$ are connected via the following formulas,
which are often used in calculations:
\begin{equation} \label{Equation 1.1}
\uptau=\log\frac{1+r}{1-r}, \quad \quad
r=\frac{e^\uptau-1}{e^\uptau+1}.
\end{equation} %

The hyperbolic area of a Borel measurable subset $E$ of $\mathbb{D}$ is %
\begin{equation} \label{Hyperbolic area}%
A_{\mathbb{D}}(E)=\int_E \lambda_{\mathbb{D}}^2(z)\,dm,
\end{equation} %
where $dm$ stands for the $2$-dimensional Lebesgue measure. In
particular,
the hyperbolic area of the disk $\mathbb{D}_r=\{z:\,|z|<r\}$,
$0<r<1$, is given by the following formula: 
$$ 
A_{\mathbb{D}}(\mathbb{D}_r)=\frac{  4 \pi \,
r^2}{1-r^2}=4\pi\sinh^2(\uptau/2).
$$ 

For the properties of geometric quantities defined above, we
recommend the monograph of A.~Beardon \cite{Be}.

Our main focus in this paper will be on the quantity, which we
call  the \emph{ conformal capacity}, or just  \emph{ capacity}
that is related to the hyperbolic capacity studied in \cite{G} and
\cite{Tsuji}.

\begin{definition} \label{Definition 1.2}
Let $E$ be a compact set in $\mathbb{D}$. The conformal
capacity 
${\rm cap}(E)$ of $E$ is defined as  %
\begin{equation}  \label{Equation 1.3}%
{\rm cap}(E)=\inf \int_{\mathbb{D}} |\nabla u|^2\,dm,
\end{equation} %
where the infimum is taken over all Lipschitz functions $u$ such
that $u=0$ on the unit circle $\mathbb{T}=\partial  \mathbb{D}$ and $u(z)\ge 1$ for $z \in E$. %
\end{definition} %

In terminology used in electrostatics, the  conformal capacity
${\rm cap}(E)$ is usually referred to as the capacity of the
condenser $(\mathbb{D},E)$ with plates $E$ and
$\mathbb{D}^*=\overline{\mathbb{C}}\setminus \mathbb{D}$ and field
$\mathbb{D}\setminus E$. Therefore, many properties and theorems
known for the capacity of a physical condenser can be applied to
the conformal capacity as well. In Section~2, we collect some of
these properties, which will be  needed for the purposes of this
work. Several available methods to prove the above mentioned
properties also will be discussed in Section~2.

\smallskip

Explicit expressions for the  conformal capacity of compact sets
are available in a few cases only. Here are three examples, in
which the conformal capacity is expressed in terms of both
Euclidean and hyperbolic characteristics of the set.

\begin{example} \label{Example 1.4}
The conformal capacity of the closure of the disk $\mathbb{D}_r$,
with Euclidean radius $r$, $0<r<1$, and hyperbolic radius
$\uptau>0,$ is given by 
$$ 
{\rm
cap}(\overline{\mathbb{D}}_r)=\frac{2\pi}{\log(1/r)}=\frac{2\pi}{\log((e^\uptau+1)/(e^\uptau-1))}.
$$ 
\end{example}

\begin{example} \label{Example 1.5}
The  conformal capacity of the interval $[-r,r]$, with hyperbolic
length equal to
$2\uptau=\ell_\mathbb{D}([-r,r])=2\log((1+r)/(1-r))$, is
$$ 
{\rm cap}([-r,r])=8\frac{\mathcal{K}(r^2)}{\mathcal{K}'(r^2)}=
8\frac{\mathcal{K}((e^\uptau-1)^2/(e^\uptau+1)^2)}{\mathcal{K}'((e^\uptau-1)^2/(e^\uptau+1)^2)}.
$$ 
Here and below, $\mathcal{K}(k)$ and
$\mathcal{K}'(k)=\mathcal{K}(\sqrt{1-k^2})$ stand for the complete
elliptic integrals of the first kind.
\end{example}

\begin{example}  \label{Example 1.7}
The conformal capacity of the interval $[0,r]$, with hyperbolic
length equal to $\uptau=\ell_\mathbb{D}([0,r])=\log((1+r)/(1-r))$,
is
\begin{equation}  \label{Equation 1.7}%
{\rm cap}([0,r])=4\frac{\mathcal{K}(r)}{\mathcal{K}'(r)}=
4\frac{\mathcal{K}((e^\uptau-1)/(e^\uptau+1))}{\mathcal{K}'((e^\uptau-1)/(e^\uptau+1))}.
\end{equation} %
\end{example} %


As Examples~\ref{Example 1.5} and \ref{Example 1.7} demonstrate,
even for compact sets as simple as a hyperbolic interval, the
conformal capacity can not be expressed in terms of elementary
functions. Thus, estimates in terms of Euclidean characteristics
of a set and numerical computations are important when working
with this capacity.

This project  originated with the following question raised by the
third-listed author of this paper. This question arose in the
course of recent work  \cite{NV1}-\cite{NRV2} and it was
experimentally studied in \cite{ENV}.

\begin{problem} \label{Problem 1.9}%
Suppose that $ 0<r<s<t<1$ and $0<u<1$ are such that the sets
$E=[0,r]\cup [s,t]$ and $E_1=[0,u]$ have equal hyperbolic lengths.
Is it true that conformal capacity of $E$ is greater than
conformal capacity of $E_1$?
\end{problem} %

It appears that this question has many interesting ramifications.
Thus, we decided to team up to discuss these questions, answer
several of them demonstrating available technique and to point out
a few remaining open questions. In the context of
Problem~\ref{Problem 1.9}, it is natural to consider compact sets
lying on any finite number of radial intervals. This is how
geometric shapes resembling animals with spikes,
and therefore the term ``hedgehog'', appeared in our study. %
  Typically, the central body $E_0$ will be a disk
$\overline{\mathbb{D}}_r$, $0<r<1$, or an emptyset and $E_k$,
$1\le k\le m$, will be a collection of closed intervals attached
to the central body. In this case, our compact sets $E$ look more
like the sea creatures called  ``stylocidaris affinis'', that is
shown in our Figure~1, than like hedgehogs as everyone knows them.
But, because the term ``hedgehog'' was already applied in the
context of our study in mathematical literature, we will stick
with it in our paper.

\begin{figure} 
\begin{minipage}{1.0\textwidth}
$$\includegraphics[scale=0.75,angle=0]{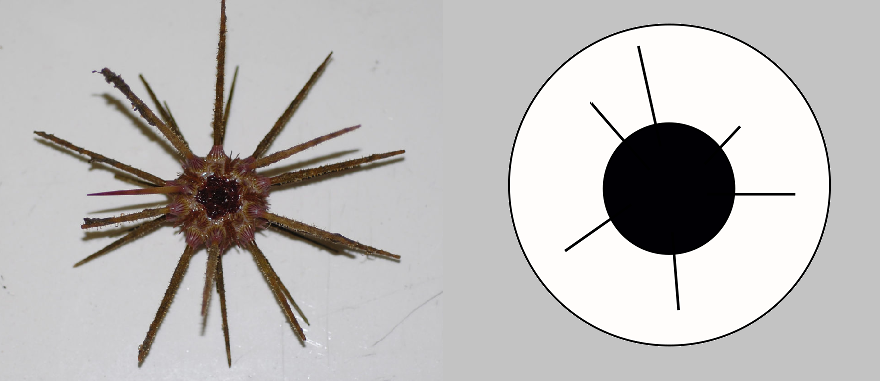}
$$
\end{minipage}
\caption{Hedgehog-like shapes: Live and geometrical.}%
\end{figure}

Our main results in Sections 3 and 4 deal with several
extremal problems for the  capacity of compact sets in the unit
disk $\mathbb{D}$, where hedgehogs possessing certain symmetry
properties play the role of the extremal configuration. Thus, in
these sections we mainly work with compact sets in $\mathbb{D}$
having components lying on a finite number of radial intervals. In
Section~3, we first demonstrate our methods on simple cases,
considered in Lemmas~\ref{Lemma 3.1} and \ref{Lemma 3.4}, when a
compact set lies on the radial segment or on the diameter of
$\mathbb{D}$. In particular, Lemma~\ref{Lemma 3.1} provides an
affirmative answer to the question stated in Problem~\ref{Problem
1.9}. Then, in several theorems presented in Section~3, we extend
our proofs to the case of compact sets lying on several radial
intervals. In Section~4, we deal with several extremal problems on
the conformal capacity for compact sets lying on a finite number
of radial intervals evenly distributed over the unit disk.


As is well known, symmetrization type transformations (such as
Steiner symmetrization, Schwarz symmetrization, P\'{o}lya circular
symmetrization, Szeg\"{o} radial symmetrization, polarization and
other) provide a standard tool to estimate capacities and many
other characteristics of sets.  Most of the classical results on
symmetrization can be found in the fundamental study by
G.~P\'{o}lya and G.~Szeg\"{o} \cite{PS}. More recent approaches to
symmetrization were developed by A.~Baernstein~II \cite{B},
V.~Dubinin  \cite{Du}, J.~Sarvas \cite{Sarvas} and also in the
papers \cite{Solynin1996a}, \cite{BrockSolynin2000},
\cite{Solynin2012} and \cite{Solynin2020}. In Section 5, we will
discuss hyperbolic counterparts of some of these transformations
and how they can be applied in problems about conformal capacity.

Finally, in Section~6, we will mention possible generalizations of
our results for conformal capacity in hyperbolic spaces of
dimension $n\ge 3$.


\section{Preliminary results on the conformal capacity} %

In this section, we recall properties of the conformal capacity
needed for our work.
We have already mentioned in the Introduction the connection of
the conformal capacity with the condenser capacity. A condenser is
a pair $(D,E)$, where $D$ is a domain in the plane and $E$ is a
compact subset of $D$. The capacity of the condenser $(D,E)$ is
defined by
\begin{equation} \label{Equation 2.1}
{\rm cap}(D,E)=\inf \int_D |\nabla u|^2\,dm,
\end{equation}
where the infimum is taken over all Lipschitz functions $u$ such
that $u\le 0$ on $\partial D$ and $u\ge 1$ on $E$.
 These functions will be called
\emph{admissible} for the condenser $(D,E)$. We want to stress
here that if $D=\mathbb{D}$,  then the infimum in (\ref{Equation
2.1}) can be taken over all admissible functions as above with an
additional requirement that $u(z)=0$ for all $z\in \mathbb{T}$.

By Theorem~3.8 of Ziemer \cite{Zim}, the capacity of the condenser
$(D,E)$ is equal to the modulus $\M(\Gamma)$ of the family
$\Gamma$ of all curves in $D\setminus E$ joining $E$ with
$\partial D$. For the definition and the basic properties of the
modulus of curve families, we refer to \cite[Chapter
II]{Jenkins1958} and \cite[Chapter 7]{HKV}.


The following invariance property of the conformal capacity, that
we often use in our proofs below, is immediate from the well-known
invariance property of the capacity of a condenser, see, for
instance, \cite[Theorem~1.12]{Du}.
\begin{proposition} \label{Proposition 2.2} %
    The conformal capacity is invariant under the M\"{o}bius self
    maps of $\mathbb{D}$ and it is invariant under reflections with
    respect to hyperbolic geodesics. Thus, if $E$ is a compact
    subset of $\mathbb{D}$ and $\varphi:\mathbb{D}\to \mathbb{D}$ is a
    M\"{o}bius automorphism or  a reflection with respect to a
    hyperbolic geodesic,
then
   $
    {\rm cap}(\varphi(E))={\rm cap}(E).
    $ 
\end{proposition} %

\medskip

Similar to the ``sets of measure zero'' in measure theory, there
are small sets that can be neglected when working with the
conformal capacity. These sets are known in the literature as
``polar sets'' or ``sets of zero logarithmic capacity''. For our
purposes the second name is more appropriate and many other
authors used it in a similar context. We recall that the
logarithmic capacity ${\rm log.cap}(E)$ of a set $E\subset
\mathbb{C}$, not necessarily
compact, is given by %
$$ %
{\rm log.cap}(E)=\sup_\mu e^{I(\mu)} \quad {\mbox{with}} \quad
I(\mu)=\iint \log|z-w|\,d\mu(z)d\mu(w),
$$
where the supremum is taken over all Borel probability measures
$\mu$ on $\mathbb{C}$  whose support is a compact subset of $E$.
For the properties of the logarithmic capacity, we refer to
Chapter~5 in T.~Ransford book \cite{R} and to the monographs
\cite{AG}, \cite{Lan}. In Proposition~\ref{Proposition 2.2.1}
below we collect results identifying ``sets of zero logarithmic
capacity'' as sets negligible for the value of the conformal
capacity. For the proofs of these results we refer to H.~Wallin's
paper \cite{Wallin}.
\begin{proposition} \label{Proposition 2.2.1} The following hold:%
\begin{enumerate} %
\item[(1)]  If $E\subset \mathbb{D}$ is compact, then ${\rm
cap}(E)=0$ if and only if ${\rm
log.cap}(E)=0$. %
\item[(2)] Let $E_1$, $E_2$ be compact sets  in $\mathbb{D}$ such
that $E_1\subset E_2$. Then ${\rm cap}(E_1)={\rm cap}(E_2)$ if and
only
if ${\rm log.cap}(E_2\setminus E_1)=0$. %
\item[(3)] If $E\subset (-1,1)$ is such that ${\rm log.cap}(E)=0$,
then $\ell_\mathbb{D}(E)=0$.
\end{enumerate} %
\end{proposition} %

We note that the set $E_2\setminus E_1$ in part (2) and the set
$E$ in part (3) of Proposition~\ref{Proposition 2.2.1} are not
necessarily compact. We stress here that the inverse statement for
part (3) of this proposition is not true, in general. For example,
if $K\subset [0,1]$ is the standard Cantor set, then its scaled
version $K_{1/3}=\{z:\,3z\in K\}$ has zero hyperbolic length but
positive logarithmic capacity, see, for example, \cite[p.143]{R}.
A compact set of logarithmic capacity zero is of zero Hausdorff
dimension \cite{Wallin}.

\medskip

Next, we will state a proposition about the existence of a
function minimizing the Dirichlet integral in
equation~(\ref{Equation 2.1}), and therefore in
equation~(\ref{Equation 1.3}) as well.
This proposition follows from classical potential theoretic
results; see, for instance,  \cite[Theorem 1]{Bag},
\cite[p.97]{Lan}.

\begin{proposition}  \label{Proposition 2.3.0} %
     Let $E$ be a compact set in $\mathbb{D}$.
     There is a unique function $u_E$, called the potential function of
     $E$, that minimizes the integral in  (\ref{Equation 2.1}); i.e. such that
    $$ 
     {\rm{cap}}(E)=\int_{\mathbb{D}} |\nabla u_E|^2\,dm.
     $$ 
Moreover, $u_E$ possesses the following properties: %
\begin{enumerate} %
\item[(1)] $u_E$ is harmonic in $\mathbb{D}\setminus E$ and
continuous on $\overline{\mathbb{D}}$,  except possibly for a
subset of $E$ of zero logarithmic capacity. %
\item[(2)] $u_E(z)=0$ for $z\in \mathbb{T}$ and $u_E(z)=1$ for all
$z\in E$,  except possibly for a subset of $E$ of zero logarithmic
capacity. %
\item[(3)] If every point of $\partial E$ is regular for the
Dirichlet problem in $\mathbb D\setminus E$, then $u_E$ is
continuous on $\overline{\mathbb D}$ and $u_E=1$ on $E$.
\end{enumerate} %
\end{proposition} %

Below are two examples of sets exceptional in the sense of parts
(1) and (2) of Proposition~\ref{Proposition 2.3.0}.

\begin{example}
Consider a compact set $E_e=\{0\}\cup\left(\cup_{n=1}^\infty
I_n\right)$, where $I_n$ is an interval $[e^{-n},e^{-n}+e^{-n^3}]$
with length $l_n=e^{-n^3}$. Since
$$ %
\sum_{n=1}^\infty \frac{n}{\log(2/l_n)}=\sum_{n=1}^\infty
\frac{n}{n^3+\log 2}<\infty,
$$ %
Wiener's criterion (see \cite[Theorem 5.4]{R}) implies that the
point $z=0\in E_e$ is irregular for the Dirichlet problem.
Therefore, the function $u_{E_e}(z)-1$ is not a barrier at $z=0$
(see \cite[Definition 4.1.4]{R}). The latter implies that the
limit $\lim_{z\to 0} u_{E_e}(z)$ does not exist. Since points
$z\in E_e$, $z\not=0$, are regular for the Dirichlet problem,
$z=0$ is the only point in $\mathbb{D}$, where $u_{E_e}$ is not
continuous.

To obtain a compact set $F_e$ such that $u_{F_e}$ has infinite
number of discontinuities, we modify our previous example as
follows. For $n\in \mathbb{N}$, let $E_e^n=\{(z+1)/2^n:\,z\in
E_e\}$. Thus, $E_e^n$ is obtained by translating and scaling the
set $E_e$. Let $F_e=\{0\}\cup \left(\cup_{n=1}^\infty
E_e^n\right)$. Our previous argument can be applied to show that
$u_{F_e}$ is not continuous at an infinite set of points
$z_n=2^{-n}$, $n\in \mathbb{N}$. Thus, $F_e$ has an infinite
subset exceptional in the sense of part (1) of
Proposition~\ref{Proposition 2.3.0}.
\end{example}

\begin{example}  %
Let $E\subset \mathbb{D}$ be a compact set of positive logarithmic
capacity, which contains a nonempty subset $E_0$, each point of
which is isolated from other points of $E$. Since $u_E$ is
harmonic and bounded, every point $z\in E_0$ is removable, which
means that $u_E$ can be extended as a function harmonic at $z$.
Since $u_E$ is not constant, it follows that $0<u_E(z)<1$ for
every point $z\in E_0$. Thus, $E_0\subset E$ is an exceptional
set, possibly infinite, as it was  mentioned in part (2) of
Proposition~\ref{Proposition 2.3.0}.
\end{example} %




\medskip

Next, we recall a subadditivity property of the conformal
capacity, which we need in the following form.

\begin{proposition} \label{Proposition 2.3.1} %
    Suppose that $E=\cup_{k=1}^n E_k$ is the union of $n\ge 2$
    compact sets $E_k$, $k=1,\ldots,n$, in $\mathbb{D}$. Then %
    \begin{equation} \label{Equation 2.8}%
       {\rm{cap}}(E)\le \sum_{k=1}^n {\rm{cap}}(E_k).
    \end{equation} %

    Moreover, if each $E_k$ has positive conformal capacity, then
(\ref{Equation 2.8}) holds with strict inequality.
\end{proposition} %

\noindent %
\emph{Proof.} %
Let $u_k$, $k=1,\ldots,n$, be an admissible function for the
condenser $(\mathbb{D},E_k)$ such that $u_k(z)=0$ for $z\in
\mathbb{T}$. Then $u(z)=\max\{u_1(z),\ldots,u_n(z)\}$ is an
admissible function for the condenser $(\mathbb{D},E)$ such that
$|u(z_1)-u(z_2)|\le \max_{1\le k\le n}|u_k(z_1)-u_k(z_2)|$ for all
$z_1,z_2\in \mathbb{D}$. The latter inequality implies that
$|\nabla u(z)|\le \max_{1\le k\le n}|\nabla u_k(z)|$ for all $z\in
\mathbb{D}$, where the gradients exist. Therefore,
\begin{equation} \label{Equation 2.8.1}
{\rm{cap}}(E)\le \int_{\mathbb{D}} |\nabla u|^2\,dm\le
\int_{\mathbb{D}} \max_{1\le k\le n}|\nabla u_k|^2\,dm\le
\sum_{k=1}^n \int_{\mathbb{D}} |\nabla u_k|^2\,dm.
\end{equation} 
Taking the infimum in this equation over all admissible functions
$u_k$, $k=1,\ldots,n$, with the properties mentioned above in this
proof, we obtain the required subadditivity property
(\ref{Equation 2.8}).

If ${\rm cap}(E_k)>0,\;k=1,2,\dots,n$, then each of the condensers
$(\mathbb{D},E_k)$ has the potential function $u_{E_k}$ that is a
non-constant harmonic function in $\mathbb D\setminus E_k$.
Therefore, $|\nabla u_{E_k}|>0$, almost everywhere in $\mathbb
D\setminus E_k$. Since $u_{E_k}$ is the potential function of
$(\mathbb{D},E_k)$ it follows that (\ref{Equation 2.8.1}) holds
with $u_k=u_{E_k}$. In this case,  the strict inequality
$\max_{1\le k\le n}|\nabla u_k(z)|<\sum_{1\le k\le n}|\nabla
u_k(z)|$ holds for all points $z$ in an annulus
$\{z:\,\rho<|z|<1\}$ with $0<\rho<1$ such that
$\{z:\,\rho<|z|<1\}\subset \mathbb{D}\setminus \cup_{1\le k\le
n}E_k$. The latter implies that if $u_k=u_{E_k}$ then the third
inequality in (\ref{Equation 2.8.1}) is strict. Therefore,
(\ref{Equation 2.8}) holds with the sign of strict inequality in
the case under consideration.\hfill $\Box$

\medskip



The proofs of our main theorems in Sections~3 and 4 rely on the
polarization technique and on the geometric interpretation of the
conformal capacity in terms of the hyperbolic transfinite
diameter.



To define the polarization of compact sets with respect to a
hyperbolic geodesic $\gamma$, we need the following terminology. %
To any hyperbolic geodesic $\gamma$, we can give an orientation by
marking one of its complementary hyperbolic halfplanes and call it
$H_+$; then the other complementary hyperbolic halfplane is given
the name $H_-$. Since every hyperbolic geodesic $\gamma$ is an arc
of a circle, we can define the classical symmetry transformation
(called also inversion or reflection) with respect to $\gamma$.
We note that the symmetry transformation with respect to a
hyperbolic geodesic $\gamma$ is a hyperbolic isometry on
$\mathbb{D}$.

\smallskip

The polarization transformation of compact sets in $\mathbb{D}$
can be defined as follows.

\begin{definition}\label{my214} %
    Let $\gamma$ be an oriented hyperbolic geodesic  in $\mathbb{D}$. Let $H_+,H_-$
    be the hyperbolic halfplanes determined by $\gamma$.
    Let $E$ be a compact set in $\mathbb{D}$ and
    let $\mathcal{R}_\gamma(E)$ denote the set symmetric to $E$ with
    respect to  $\gamma$. The polarization $\mathcal{P}_\gamma(E)$ of $E$ with
    respect to $\gamma$ is defined by
    \begin{equation}  \label{Equation 2.10}%
    \mathcal{P}_\gamma(E)=((E\cup \mathcal{R}_\gamma(E))\cap \overline{H_+})\cup((E\cap \mathcal{R}_\gamma(E))\cap H_-).
    \end{equation}
\end{definition} %

Equation (\ref{Equation 2.10}) can be also written in the form
$$ 
\mathcal{P}_\gamma(E)=((E\cup \mathcal{R}_\gamma(E))\setminus
H_-)\cup ((E\cap \mathcal{R}_\gamma(E))\setminus H_+).
$$ 

The polarization transformation was introduced by V.~Wolontis in
1952, \cite{Wolontis}. Wolontis' work remained unnoticed until
1984, when V.~N.~Dubinin used this transformation to solve
A.~A.~Gonchar's problem on the capacity of a condenser with plates
on a fixed straight line interval. The name ``polarization'' was
also  suggested by Dubinin \cite{Du85}. The following proposition,
describing the change of the conformal capacity under
polarization, is an important ingredient of the proofs in the
following sections, see \cite[Theorem 3.4]{Du}, \cite[Theorem
2.8]{BP2012}.

\begin{proposition} \label{Proposition 2.5} %
    Let $E$ be a compact set in $\mathbb{D}$ and $\mathcal{P}_\gamma(E)$ be the
    polarization of $E$ with respect to an oriented hyperbolic geodesic $\gamma$.  Then  %
    \begin{equation} \label{Equation 2.12}%
    {\rm cap}(\mathcal{P}_\gamma(E))\le
    {\rm cap}(E).
    \end{equation} %

    Furthermore, equality occurs in (\ref{Equation 2.12}) if and
    only if $\mathcal{P}_\gamma(E)$ coincides with $E$ up to reflection with respect to $\gamma$
    and up to a set of zero logarithmic capacity. %
\end{proposition} %

When working with hedgehog structures, the following particular
case of Proposition~\ref{Proposition 2.5} is useful.

\begin{corollary} \label{Corollary 2.1} %
Under the assumptions of Proposition~\ref{Proposition 2.5}, let
$E$ be the closure of the union of a finite or infinite number of
non-overlapping closed intervals on the diameter $(-1,1)$. Then
(\ref{Equation 2.12}) holds with the sign of strict inequality
unless $\mathcal{P}_\gamma(E)$ coincides with $E$ up to reflection
with respect to $\gamma$.
\end{corollary} %


\medskip



 One more useful characteristic of compact sets in the hyperbolic plane, the hyperbolic transfinite diameter,
 was introduced by M.~Tsuji \cite{Tsuji}. It is defined as
follows (see \cite[p.94]{Tsuji} or \cite[Section 1.4]{Du}).

\begin{definition}\label{mydef} %
    Let $E\subset \mathbb{D}$ be a compact set in $\mathbb{D}$. The hyperbolic transfinite diameter of $E$ is defined as
    \begin{equation}  \label{Hyperbolic transfinite diameter}%
    {\rm d}_h(E)=\lim_{n\to \infty} \max \prod_{1\leq j<k\leq n} \left(
    p_{\mathbb D}(z_j,z_k)\right)^{2/[n(n-1)]},
   \end{equation}
    where $p_\mathbb{D}(z_j,z_k)$ stands for the pseudo-hyperbolic metric defined by (\ref{Pseudo-hyperbolic distance})
    and the maximum is taken of all $n$-tuples of points $z_1,\ldots,z_n$ in $E$. %
\end{definition} %

The following relation was established in 
\cite{Tsuji}. %

\begin{proposition} \label{Proposition 2.20.0} %
    Let $E$ be a compact set in $\mathbb{D}$. Then
    \begin{equation}  \label{Capacity via transfinite diameter}
    {\rm cap}(E)=\left [-\frac{1}{2\pi}\log{\rm d}_h(E)\right ]^{-1}.
   \end{equation} 
\end{proposition} %

Let $E$ be a compact set in $\mathbb{D}$. A map $\varphi:E\to
{\mathbb D}$ is called a hyperbolic contraction on $E$ if for
every $z_1,z_2\in E$,
$$
d_{\mathbb D}(\varphi(z_1),\varphi(z_2))\leq d_{\mathbb
D}(z_1,z_2),
$$
where $d_\mathbb{D}(\cdot,\cdot)$ stands for the hyperbolic metric
defined by (\ref{Hyperbolic distance}). Furthermore, $\varphi:E\to
{\mathbb D}$ is called a strict hyperbolic contraction on $E$ if
there is $k$, $0<k<1$, such that for every $z_1,z_2\in E$,
$$
d_{\mathbb D}(\varphi(z_1),\varphi(z_2))\leq k\,d_{\mathbb
D}(z_1,z_2).
$$

The following contraction principle is immediate from the
Definition~\ref{mydef} and Proposition~\ref{Proposition 2.20.0}.

\begin{proposition} \label{Proposition 2.19} %
    Let $E$ be a compact set in $\mathbb{D}$. Let $\varphi:E\to
    \mathbb{D}$ be a hyperbolic contraction.  Then ${\rm
        cap}(\varphi(E))\le {\rm cap}(E)$. %

        Moreover, if $\varphi$ is a strict hyperbolic contraction on $E$,
        then ${\rm
        cap}(\varphi(E))< {\rm cap}(E)$. %
\end{proposition} %



\medskip

\begin{rem} %
We note here that the polarization transformation is not
contracting, in general. For example, polarizing the set $E=\{\pm
i/4,\pm (1-i)/4\}\subset \mathbb{D}$ with respect to the diameter
$I=(-1,1)$ with its standard orientation, we obtain the polarized
set ${\mathcal P}_I(E)=\{\pm i/4,(\pm 1+i)/4\}$. Then for every
one-to-one map $\varphi:E\to {\mathcal P}_I(E)$, there is a pair
of points $z_1,z_2\in E$ such that $d_{\mathbb
D}(\varphi(z_2),\varphi(z_1))>d_{\mathbb D}(z_2,z_1)$; one can
easily verify this inequality by considering each one of the
possible maps $\varphi$.
\end{rem} %



To study the limit behavior of the conformal capacity
${\rm{cap}}(E)$, when some of the components of $E$ tend to the
boundary of $\mathbb{D}$, we need a hyperbolic analog of the
dispersion property of the Newtonian capacity discussed in
\cite{Solynin2021}. Let $E_1,\ldots,E_n$ be disjoint nonempty
compact sets
in $\mathbb{D}$, 
not
necessarily connected, and let $E=\cup_{k=1}^n E_k$.

\begin{definition} \label{Definition 2.1}%
    By a hyperbolic dispersion of $E=\cup_{k=1}^n E_k$ we mean a mapping
    $\varphi:E\times [0,\infty)\to\mathbb{D}$ satisfying the following
    properties:
    \begin{enumerate} %
        \item[(1)] For each $k$, the restriction
        $\varphi:E_k\times[0,\infty)\to\mathbb{D}$ is a rigid hyperbolic motion of
        $E_k$, which depends continuously on the parameter $t\in
        [0,\infty)$, such that
        $\varphi(x,0)=x$ for all $x\in E$. %
        \item[(2)] If $0\le t_1<t_2$, then for each $k$ and $j$, $k\not =
        j$, the hyperbolic distances between the images $\varphi(E_k,t)$ and
        $\varphi(E_j,t)$ satisfy the following inequalities:%
        $$ %
        d_\mathbb{D}\,(\varphi(E_k,t_1),\varphi(E_j,t_1))\le
        d_\mathbb{D}\,(\varphi(E_k,t_2),\varphi(E_j,t_2)).
        $$ %
        \item[(3)] For each $k$ and $j$, $k\not =
        j$, %
        $$ %
        d_\mathbb{D}\,(\varphi(E_k,t),\varphi(E_j,t))\to \infty \quad
        {\mbox{as $t\to \infty$.}}
        $$ %
    \end{enumerate} %
\end{definition} %

Thus, hyperbolic  dispersion of $E$ is a process moving the
subsets $E_1,\ldots,E_n$  farther and farther from each other,
resembling the scattering of galaxies of our Universe.

We stress here that not every finite collection of compact sets
admits hyperbolic dispersion. For example, the set $E=E_1\cup
E_2$, where $E_1=\{0\}$ and $E_2=\{z=re^{i\theta}:\,|\theta|\le
\pi-\varepsilon\}$ with $0<r<1$ and sufficiently small
$\varepsilon>0$, cannot be hyperbolically dispersed in the sense
of Definition~\ref{Definition 2.1}. On the other hand, any union
$E=E_1\cup E_2$ of two non-intersecting compact sets, each of
which lies on a radial interval, can be hyperbolically dispersed.

\medskip

The following useful result is a hyperbolic counterpart of
Proposition~5 proved in \cite{Solynin2021}.

\begin{proposition}   \label{Proposition 2.3.2} %
    Let $\varphi:E\times [0,\infty)\to\mathbb{D}$ be a hyperbolic dispersion of a
    compact
    set $E=\cup_{k=1}^n E_k$, as above.   Then %
    \begin{equation} \label{Equation 2.20}%
    {\rm{cap}}(\varphi(E,t))\to  \sum_{k=1}^n {\rm{cap}}(E_k), \quad {\mbox{as $t\to \infty$.}}
    \end{equation} %
\end{proposition} %

\smallskip

In the proof of  Proposition~\ref{Proposition 2.3.2}, we will need
the following elementary arithmetic result.

\begin{lemma} \label{Arithmetic approximation} %
Let $0<\alpha_k<1$, $k=1,\ldots,n$, be such that
$\sum_{k=1}^n\alpha_k=1$. Then there are $k$ sequences of positive
integers $m_{j,k}$, $j=1,2,\ldots$, such that if $m_j=\sum_{k=1}^n
m_{j,k}$, then $m_j\to \infty$ and $m_{j,k}/m_j\to \alpha_k$ as
$j\to \infty$.
\end{lemma} %

\noindent %
\emph{Proof.} 
Consider rational approximations of
$\alpha_k$, $k=1,\ldots,n-1$; i.e. consider $(n-1)$ sequences %
$$  %
\frac{a_{j,k}}{b_{j,k}}\to \alpha_k, \quad \ k=1,\ldots,n-1, \
{\mbox{ as $j\to \infty$,}}
$$ %
where $a_{j,k}$, $b_{j,k}$ are positive integers. Then consider
the
sequence %
$$ %
m_j=j\prod_{k=1}^{n-1}b_{j,k}\to \infty
$$ %
and the sequences %
$$ %
m_{j,k}=\frac{a_{j,k}m_j}{b_{j,k}}, \quad k=1,\ldots,n-1.
$$ %
Clearly, %
$$ %
\frac{m_{j,k}}{m_j}=\frac{a_{j,k}}{b_{j,k}}\to \alpha_k \quad
{\mbox{as $j\to \infty$.}}
$$ %
Also, we put %
$$ %
m_{j,n}=m_j-\sum_{k=1}^{n-1}m_{j,k}.
$$ %
Then %
$$ %
\frac{m_{j,n}}{m_j}=1-\sum_{k=1}^{n-1}\frac{m_{j,k}}{m_j}\to
1-\sum_{k=1}^{n-1}\alpha_k=\alpha_n>0 \quad {\mbox{as $j\to
\infty$.}}
$$ %
The latter relation shows that $m_{j,n}>0$ for all $j$
sufficiently large. Therefore, we can remove a finite number of
terms from the sequence $m_j$ and from the sequences $m_{j,k}$ and
then re-enumerate these sequences to obtain sequences with the
required properties. \hfill $\Box$

\medskip

\noindent %
\emph{Proof of Proposition~\ref{Proposition 2.3.2}.} For
$k=1,2,\dots, n$ and $t\geq 0$, we set $E^t=\varphi(E,t)$,
$E_k^t=\varphi(E_k,t)$. Since the conformal
capacity is invariant under hyperbolic motions, %
${\rm{cap}}(E_k^t)={\rm{cap}}(E_k)$
for all $k=1,\ldots,n$ and all $t\ge 0$. This together with the
subadditivity property of Proposition~\ref{Proposition 2.19}
implies that
\begin{equation}\label{Equation 2.001}
{\rm cap}(E^t)\leq \sum_{k=1}^n {\rm
    cap}(E_k^t)=\sum_{k=1}^n {\rm cap}(E_k).
\end{equation}

We assume without loss of generality  that ${\rm{cap}}(E_k)>0$ for
all $k=1,\ldots,n$. Then we set
$\alpha_k={\rm{cap}}(E_k)/\sum_{k=1}^n {\rm{cap}}(E_k)$ and will
use the sequences $m_{j,k}$, $k=1,\ldots,n$, and $m_j$ defined as
in the proof of Lemma~\ref{Arithmetic approximation} for our
choice of $\alpha_k$.

Let $z_{j,k}^s$, $s=1,\ldots,m_{j,k}$, be points in $E_k$ such
that %
\begin{equation}\label{Equation 2.002}
\prod_{1\le l<s\le m_{j,k}}p_\mathbb{D}(z_{j,k}^l,z_{j,k}^s)=\max
\prod_{1\le l<s\le m_{j,k}}p_\mathbb{D}(z^l,z^s),
\end{equation}
where the maximum is taken over all $m_{j,k}$-tuples of points
$z^1,z^2,\ldots,z^{m_{j,k}}$ in $E_k$. For $k=1,\ldots,n$,
$j=1,2,\ldots$, $s=1,\ldots,m_{j,k}$, and $t\ge 0$, we set
$z_{j,k}^{s,t}=\varphi(z_{j,k}^s,t)$. Since $\varphi$ is a
hyperbolic motion on each $E_k$, we have %
\begin{equation}\label{Equation 2.003}
p_\mathbb{D}(z_{j,k}^{l,t},z_{j,k}^{s,t})=p_\mathbb{D}(z_{j,k}^l,z_{j,k}^s),
\end{equation} %
for all points  $z_{j,k}^l$, $z_{j,k}^s$ defined above and all
$t\ge 0$.

For our choice of points, it follows from (\ref{Equation 2.003})
and equation~(\ref{Hyperbolic transfinite diameter})
of Definition~\ref{mydef} that %
\begin{equation}\label{Equation 2.004}
d_h(E^t)\ge \limsup_{j\to \infty} \left[\Pi_j^t\, \prod_{k=1}^n
\Pi_{j,k}\right]^{2/m_j(m_j-1)},
\end{equation}
where %
\begin{equation}\label{Equation 2.005}
\Pi_{j,k}=\prod_{1\le l<s\le m_{j,k}}
p_\mathbb{D}(z_{j,k}^l,z_{j,k}^s),
\end{equation} %
and %
\begin{equation}\label{Equation 2.006}
\Pi_j^t=\prod p_\mathbb{D}(z_{j,k_1}^{l,t},z_{j,k_2}^{s,t}),
\end{equation} %
where the product in (\ref{Equation 2.006}) is taken over all
pairs of points $z_{j,k_1}^{l,t}$, $z_{j,k_2}^{s,t}$ such that
$1\le l\le m_{j,k_1}$, $1\le s\le m_{j,k_2}$ and $k_1\not=k_2$.

Using equations~(\ref{Hyperbolic transfinite diameter}),
(\ref{Equation 2.002}), and (\ref{Equation 2.005}) and taking into
account our choice of points $z_{j,k}^s$, $1\le s\le m_{j,k}$ and
the limit relation
$\lim_{j\to \infty} m_{j,k}/m_j= \alpha_k$, we conclude that %
\begin{align}\label{Equation 2.007}
\lim_{j\to
\infty}\left(\Pi_{j,k}\right)^{2/m_j(m_j-1)}&=\lim_{j\to
\infty}\left(\prod_{1\le l<s\le m_{j,k}}
p_\mathbb{D}(z_{j,k}^l,z_{j,k}^s)\right)^{\frac{2}{m_{j,k}(m_{j,k}-1)}\cdot\frac{m_{j,k}(m_{j,k}-1)}{m_j(m_j-1)}}\\
&=\left(d_h(E_k)\right)^{\alpha_k^2}. \nonumber
\end{align} %

Our assumption that $d_\mathbb{D}(E_{k_1}^t,E_{k_2}^t)\to \infty$
when $k_1\not=k_2$ and $t\to \infty$ and relations
(\ref{Hyperbolic distance}), (\ref{Pseudo-hyperbolic distance}),
imply that for every $\varepsilon>0$ there exists $t_\varepsilon
>0$ such that if $k_1\not=k_2$ then %
$$ %
p_\mathbb{D}(z_{j,k_1}^{l,t},z_{j,k_2}^{s,t})>1-\varepsilon \quad
{\mbox{for all $t\ge t_\varepsilon$.}}
$$ %
This inequality together with (\ref{Equation 2.006}) imply that
\begin{equation}\label{Equation 2.008}
\left(\Pi_j^t\right)^{2/m_j(m_j-1)}\ge 1-\varepsilon \quad
{\mbox{for all $t\ge t_\varepsilon$.}}
\end{equation} %

Combining (\ref{Equation 2.004}), (\ref{Equation
2.007}), and (\ref{Equation 2.008}), we obtain the following: %
$$  %
d_h(E^t)\ge (1-\varepsilon)\, \prod_{k=1}^n
\left(d_h(E_k)\right)^{\alpha_k^2} \quad {\mbox{for all $t\ge
t_\varepsilon$.}}
$$ %
The latter inequality together with (\ref{Capacity via transfinite
diameter}) implies that for all $t\ge t_\varepsilon$, %
\begin{align}\label{Equation 2.009}
{\rm{cap}}(E^t)&=\frac{1}{-\frac{1}{2\pi}\log d_h(E^t)}\ge
\frac{1}{-\frac{1}{2\pi}\log \left((1-\varepsilon)\prod_{k=1}^n
(d_h(E_k))^{\alpha_k^2}\right)}\\
&=\frac{1}{\sum_{k=1}^n
(\alpha_k^2/{\rm{cap}}(E_k))-\frac{1}{2\pi}\log (1-\varepsilon)}
\nonumber \\
&=\frac{\sum_{k=1}^n
{\rm{cap}}(E_k)}{1-\frac{1}{2\pi}\log(1-\varepsilon)\sum_{k=1}^n
{\rm{cap}}(E_k)}. \nonumber
\end{align} %
Since $\varepsilon>0$ can be chosen arbitrarily small, it follows
from (\ref{Equation 2.009}) that %
\begin{equation}\label{Equation 2.0010}
\liminf_{t\to \infty}{\rm{cap}}(E^t)\ge \sum_{k=1}^n
{\rm{cap}}(E_k).
\end{equation} %
Finally, equations (\ref{Equation 2.001}) and (\ref{Equation
2.0010}) imply (\ref{Equation 2.20}). \hfill  $\Box$

\medskip


\begin{rem}  %
We note here that the conformal capacity of a compact set
$E\subset \mathbb{D}$ is not monotone under hyperbolic dispersion,
in general.
\end{rem}

To give an example of such non-monotonicity, we consider a family
of hedgehogs $E(t)$, $t\ge 0$, with $E(t)$ consisting of a fixed
central body $C_{r_0}(\alpha)=\{r_0e^{i\theta}:\,|\theta|\le
\alpha\}$, $0<r_0<1$, $0<\alpha<\pi$, and single varying spike
$E(t)=[r_1(t)e^{i(\alpha-t)},r_2(t)e^{i(\alpha-t)}]$, which we
define as follows.

We put $r_1(0)=r_1$, $r_2(0)=r_2$, where $r_0<r_1<r_2<1$. Using
polarization with respect to appropriate hyperbolic geodesics
$\gamma=(-e^{i\theta},e^{i\theta})$, we find that ${\rm
cap}(C_{r_0}(\alpha)\cup [r_1e^{i(\alpha-t)},r_2e^{i(\alpha-t)}])$
strictly decreases, when $t$ varies from $0$ to $\alpha$. At the
same time
$d_\mathbb{D}(C_{r_0}(\alpha),[r_1e^{i(\alpha-t)},r_2e^{i(\alpha-t)}])$
and $\ell_\mathbb{D}([r_1e^{i(\alpha-t)},r_2e^{i(\alpha-t)}])$ are
constant for $0\le t\le \alpha$.  Using these properties and the
well-known convergence result, which is stated in
Proposition~\ref{Proposition 2.31} below, we conclude that there
is a strictly increasing function $r_1(t)$ such that $r_1(0)=r_1$,
$r_1<r_1(\alpha)<1$, and a function $r_2(t)$, $r_1(t)<r_2(t)<1$,
such that ${\rm cap}(C_{r_0}(\alpha)\cup
[r_1(t)e^{i(\alpha-t)},r_2(t)e^{i(\alpha-t)}])$ strictly
decreases,
$d_\mathbb{D}(C_{r_0}(\alpha),[r_1(t)e^{i(\alpha-t)},r_2(t)e^{i(\alpha-t)}])$
strictly increases, while the hyperbolic length
$\ell_\mathbb{D}([r_1(t)e^{i(\alpha-t)},r_2(t)e^{i(\alpha-t)}])$
remains constant on $0\le t\le \alpha$.

Now, we put $E(t)=C_{r_0}(\alpha)\cup
[r_1(t)e^{i(\alpha-t)},r_2(t)e^{i(\alpha-t)}]$ for $0\le t\le
\alpha$ and, for $t\ge \alpha$, we define $E(t)$ as
$C_{r_0}(\alpha)\cup [r_1(t),r_2(t)]$ with
$r_1(t)=(t-\alpha+\alpha r_1(\alpha)/t$ and $r_2(t)$ such that
$\ell_\mathbb{D}([r_1(t),r_2(t)])=\ell_\mathbb{D}([r_1,r_2])$. The
family of hedgehogs $E(t)$ defines a dispersion of compact sets
$C_{r_0}(\alpha)$ and $[r_1e^{i\alpha},r_2e^{i\alpha}]$ such that
${\rm cap}(E(t))$ strictly decreases on $0\le t\le \alpha$.
Furthermore, using polarization with respect to appropriate
geodesics, as we will demonstrate it later in the proof of
Lemma~\ref{Lemma 3.4}, one can show that ${\rm cap}(E(t))$
strictly increases on the interval $t\ge \alpha$.

\medskip

For our proofs, we need two results on the sequences of compact
sets in $\mathbb{D}$ convergent in an appropriate sense.

\begin{proposition}[see, {\cite[Theorem 1.11]{Du}}]  \label{Proposition 2.31} %
Let $E_k$, $k=1,2,\ldots$, be a sequence of compact sets in
$\mathbb{D}$, such that $E_{k+1}\subset E_k$ for all
$k=1,2,\ldots$, and
let $E=\cap_{k=1}^\infty E_k$. Then %
$$ 
{\rm cap}(E_k) \to {\rm cap}(E) \quad {\mbox{as $k\to \infty$.}}
$$ 
\end{proposition} %

To state our next proposition, we recall that the Hausdorff
distance between two compact sets $K,L$ in the plane is given by
$$ d_{\rm H}(K,L)=\max\{{\rm dist}(x,K),{\rm dist}(y,L): x\in
L,\;y\in K\}.
$$
The following convergence result follows from \cite[Theorem
7]{Ase}.

\begin{proposition} \label{Proposition 2.33} %
For fixed  $\delta>0$, let $E_k$, $k=1,2,\ldots$, be a sequence of
compact sets on the diameter $(-1,1)$, each of which consists of a
finite number of closed intervals such that the hyperbolic length
of each of these intervals is $\ge \delta$. If the sequence $E_k$
converges in the Hausdorff metric to a compact
set $E\subset \mathbb{D}$, then %
$$ 
{\rm cap}(E_k) \to {\rm cap}(E) \quad {\mbox{as $k\to \infty$.}}
$$ 
\end{proposition} %

\section{Hedgehogs with geometric restrictions on the number of spikes} 

We start with the following monotonicity result, which, in
particular, answers the question raised in Problem~\ref{Problem
1.9}.

\begin{lemma} \label{Lemma 3.1}
    Suppose that  $-1<a<1$ and $\uptau>0$ are
    fixed and $b$ varies in the interval $[a,1)$. Let $c=c(b)$,
    $b<c<1$, be such that $\ell_{\mathbb{D}}([b,c])=\uptau$. Let $E_0\subset (-1,a]$
    be a compact set consisting of a finite number of non-degenerate intervals
and let
    $E(b)=E_0\cup [b,c(b)]$. Then ${\rm cap}(E(b))$ is a continuous function that strictly
    increases from ${\rm
        cap}(E_0\cup[a,c(a)])$ to ${\rm cap}(E_0)+{\rm cap}([a,c(a)])$, when $b$ runs from $a$ to $1$. %
\end{lemma} %

\noindent %
\emph{Proof.} The continuity property of ${\rm cap}(E(b))$ follows
from Proposition~\ref{Proposition 2.33}. To prove the monotonicity
of ${\rm cap}(E(b))$, we consider $b_1$, $b_2$ such that $a\leq
b_1<b_2<1$ and note that $c(b_1)<c(b_2)$. Let $\gamma$ be a
hyperbolic geodesic that is orthogonal to the hyperbolic interval
$[b_1,c(b_2)]_h$ at its midpoint.  We give an orientation to
$\gamma$ by marking its complementary hyperbolic halfplane $H_+$
with $b_1\in H_+$, see Figure~2, which illustrates the proof of
this lemma. Notice that under our assumptions, $E_0\subset H_+$
and, since reflections with respect to hyperbolic geodesics
preserve hyperbolic lengthes,  the hyperbolic interval
$I_1=[b_1,c(b_1)]_h$ coincides with the reflection of
$I_2=[b_2,c(b_2)]_h$ with respect to $\gamma$. Therefore, the
polarization $\mathcal{P}_\gamma(E(b_2))$ of $E(b_2)$ with respect
to $\gamma$ coincides with the set $E(b_1)$ if $b_1\not\in E_0$
and with the set $E(b_1)\cup\{c(b_2)\}$ otherwise, and the set
$\mathcal{P}_\gamma(E(b_2))\setminus E(b_2)=(c(b_1),c(b_2))_h$ is
a non-degenerate interval and thus it has positive logarithmic
capacity. Furthermore, since $E_0\subset H_+$, it follows that the
set $\mathcal{P}_\gamma(E(b_2))$ differs from the reflection of
$E(b_2)$ with respect to $\gamma$ by a set of positive logarithmic
capacity.
 So, applying Proposition
\ref{Proposition 2.5}, we conclude that ${\rm
    cap}(E(b_1))<{\rm cap}(E(b_2))$. Thus we proved that the function
${\rm cap}(E(b))$ is  strictly increasing. The assertion about the
range of this function follows from the dispersion property of
Proposition \ref{Proposition 2.3.2} and from the convergence
property
 stated in Proposition~\ref{Proposition 2.33}. \hfill $\Box$

\begin{rem} %
The proof of Lemma~\ref{Lemma 3.1} remains valid if $E_0$ is any
compact set in the hyperbolic halfplane $H_+$ defined as in the
proof above for the hyperbolic geodesic $\gamma$ passing through
the point $a$ such that the set $E(b)=E_0\cup [b,c(b)]_h$
satisfies the assumptions of Proposition~\ref{Proposition 2.33}.
\end{rem} %

\begin{figure} 
\begin{minipage}{1.0\textwidth}
$$\includegraphics[scale=0.25,angle=0]{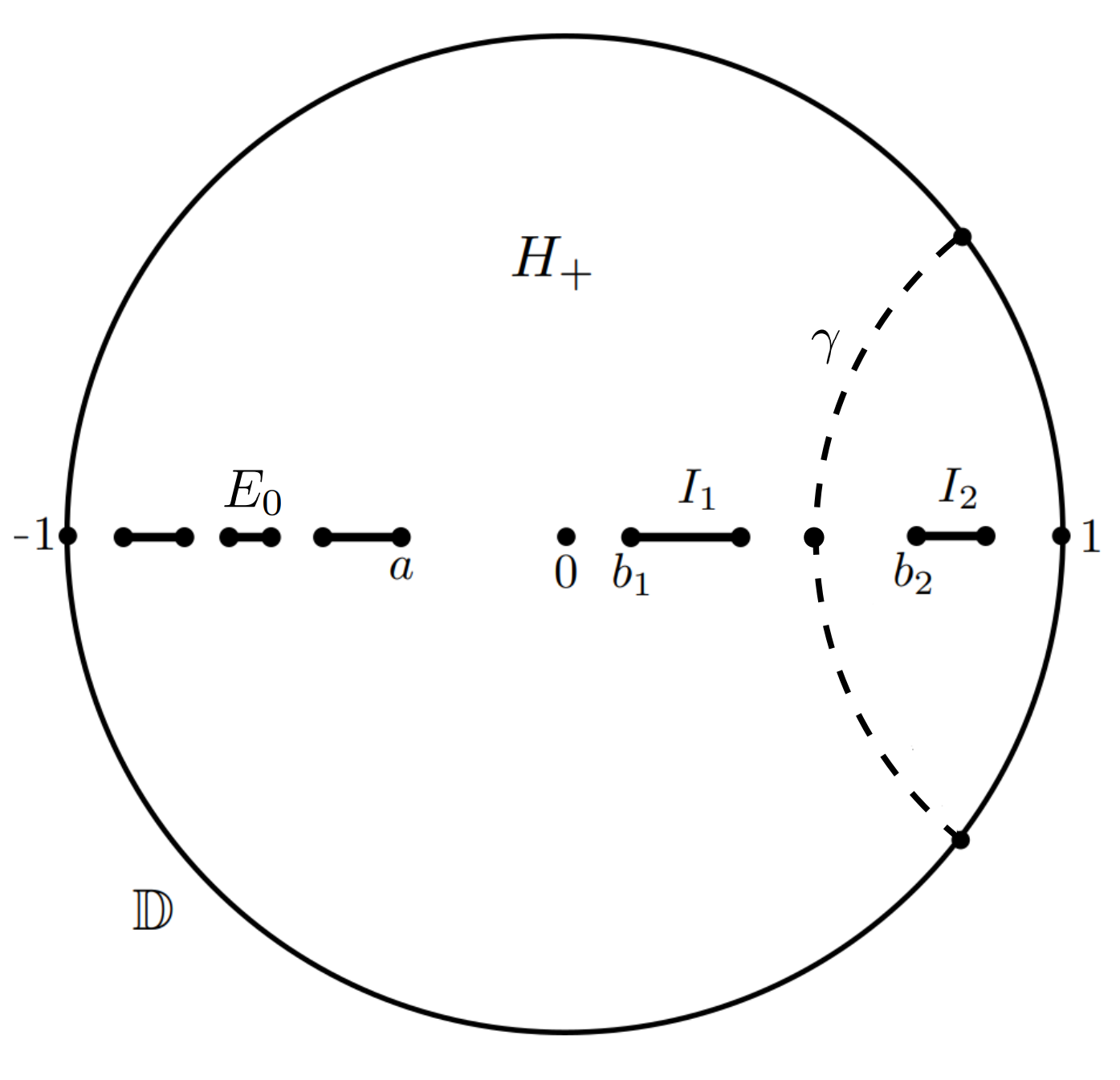}
$$
\end{minipage}
\caption{Hedgehog with one moving interval.}%
\end{figure}

\begin{rem} %
The non-strict monotonicity property in Lemma~\ref{Lemma 3.1} also
follows from the contraction principle of Proposition~\ref{Proposition 2.19}. %
\end{rem} %

 \medskip

 Next, we will use Lemma~\ref{Lemma 3.1} to prove a lower bound
 for the conformal capacity of compact sets lying on the diameter
 of $\mathbb{D}$.

 \begin{lemma}\label{Lemma 3.4}
 Let $\uptau>0$ and let $r(\uptau)$ be defined as in (\ref{Equation 1.1}).
  If $E\subset (-1,1)$ is a compact set such
that $\ell_{\mathbb{D}}(E)=\uptau$, then %
\begin{equation} \label{Equation 3.5}
{\rm cap}([0,r(\uptau)])\le {\rm cap}(E). %
\end{equation} %
Equality occurs here if and only if $E$ coincides with some
interval
$[a,b]\subset (-1,1)$ up to a set of zero logarithmic capacity.  %
\end{lemma} %

\noindent %
\emph{Proof.} (a) Suppose first that $E$ consists of $n\ge 2$
intervals $[a_k,b_k]$, $0=a_1<b_1<a_2<b_2<\cdots<a_n<b_n<1$. Let
$E_1$ be the compact set obtained from $E$ by replacing the pair
of intervals $[a_{n-1},b_{n-1}]$ and $[a_n,b_n]$ with a single
varying interval $[a_{n-1},b'_{n-1}]$ such that
$\ell_\mathbb{D}([a_{n-1},b'_{n-1}])=\ell_\mathbb{D}([a_{n-1},b_{n-1}])
+\ell_\mathbb{D}([a_n,b_n])$. It follows from the monotonicity
property of Lemma~\ref{Lemma 3.1} that ${\rm cap}(E_1)<{\rm
cap}(E)$. Applying this procedure of merging two intervals into a
single interval $n-1$ times, we obtain the
inequality~(\ref{Equation 3.5}) with the sign of strict
inequality.

(b) If $E$ is a more general compact set, not the union of a
finite number of intervals, we proceed as follows. Since the
subset of isolated points of $E$ has zero logarithmic capacity we
can remove it without changing the conformal capacity and the
hyperbolic length of $E$. Thus, we assume that $E$ does not have
isolated points. Also, since both the conformal capacity and
hyperbolic length are invariant under conformal automorphisms of
$\mathbb{D}$, we may assume that $\min\{\Re z:\,z\in E\}=0$.

We will use the following approximation argument. The set
$(-1,1)\setminus E$ is an open subset of $(-1,1)$ and therefore,
in the case under consideration,  it is a countably infinite union
of open disjoint intervals $I_k$, $k=1,2,\ldots$. We enumerate
these intervals such that $I_1=(-1,0)$ and $I_2$ has one of its
end points at $1$. Setting $E_n=(-1,1)\setminus \cup_{k=1}^{n+1}
I_k$, $n=1,2,\dots$, we obtain a sequence of compact sets $E_n$,
each consisting of a finite number of disjoint closed intervals on
$(-1,1)$, such that $E_{n+1}\subset E_n$ for all $n$ and
$E=\cap_{n=1}^\infty E_n$. So $E$ is approximated by a finite
union of closed intervals. Therefore, $\lim_{n\to \infty}
\ell_\mathbb{D}(E_n)\to \ell_\mathbb{D}(E)$ and, by
Proposition~\ref{Proposition 2.31}, $\lim_{n\to\infty}{\rm
cap}(E_n)={\rm cap}(E)$.

Let $F_n=[0,a_n]$ be a closed interval on $[0,1)$ such that
$\ell_\mathbb{D}(F_n)=\ell_\mathbb{D}(E_n)$. Then, by
part (a) of this proof, %
\begin{equation}   \label{Equation 3.6} %
{\rm cap}(F_n)<{\rm cap}(E_n). %
\end{equation} %
Furthermore, $F_{n+1}\subset F_n$ for all $n$ and
$\cap_{n=1}^\infty F_n=[0,r(\uptau)]$. Thus, ${\rm cap}(F_n)\to
{\rm cap}([0,r(\uptau)])$, by Proposition~\ref{Proposition 2.31}.
Therefore, passing to the limit in (\ref{Equation 3.6}), we obtain
(\ref{Equation 3.5}).

(c) Here we prove the equality statement. If $E$ coincides with
some interval $[a,b]$ up to a set of zero logarithmic capacity
then, by Proposition~\ref{Proposition 2.2.1},
$\ell_\mathbb{D}([a,b])=\ell_\mathbb{D}(E)=\uptau$ and ${\rm
cap}([a,b])={\rm cap}(E)$. This together with the conformal
invariance property of the capacity implies that ${\rm
cap}(E)={\rm cap}([0,r(\uptau)])$.

Suppose now that (\ref{Equation 3.5}) holds with the sign of
equality. Let $a$ and $b$ denote the infimum and supremum of the
set of Lebesgue density points of $E$. We may assume, without loss
of generality, that $a=0$.

If $r(\uptau)<b$,  then the set $[0,r(\uptau)]\setminus E$
contains an open interval $(c_1,c_2)$ with $0<c_1<c_2<r(\uptau)$.
Let $\mathcal{P}_\gamma(E)$ denote the polarization of $E$ with
respect to a hyperbolic geodesic $\gamma$ that is orthogonal to
the hyperbolic interval $[(c_1+c_2)/2,b]_h$ at its midpoint and
oriented such that $0\in H_+$. Notice, that under our assumptions,
the set $\mathcal{P}_\gamma(E)$ differs from $E$ by a set of
positive one-dimensional Lebesgue measure, and therefore by a set
of positive logarithmic capacity, and also it differs from the
reflection of $E$ with respect to $\gamma$ by a set of positive
one-dimensional Lebesgue measure, and therefore by a set of
positive logarithmic capacity. Hence, by
Proposition~\ref{Proposition 2.5},
\begin{equation}  \label{Equation 3.7}%
{\rm cap}(\mathcal{P}_\gamma(E))<{\rm cap}(E). %
\end{equation} %
 Since the polarization with respect to hyperbolic geodesics
preserves hyperbolic length, we have
$l_\mathbb{D}(\mathcal{P}_\gamma(E))=\uptau$. Therefore, it
follows from the assumption ${\rm cap}(E)={\rm
cap}([0,r(\uptau)])$ and our proof in parts (a) and (b) above,
that ${\rm cap}(E)\le {\rm cap}(P_\gamma(E)$, which contradicts
equation~(\ref{Equation 3.7}). Since the assumption that
$r(\uptau)<b$ leads to a contradiction, we must have
$b=r(\uptau)$.

In the latter case, $[0,r(\uptau)]\subset E$ and, since ${\rm
cap}(E)={\rm cap}([0,r(\uptau)])$, it follows from
Proposition~\ref{Proposition 2.2.1} that $E\setminus
[0,r(\uptau)]$ has zero logarithmic capacity, which completes the
proof of the lemma. \hfill $\Box$

\medskip

Actually, our proof of Lemma~\ref{Lemma 3.4} gives us a more
general result, which we state as the following corollary.

\begin{corollary}\label{Corollary 3.1}
For $\uptau>0$ and $-1<a<1$, let $\rho(a,\uptau)\in(a,1)$ be such
that $\ell_\mathbb{D}([a,\rho(a,\uptau)])=\uptau$. Let $\gamma$ be
a hyperbolic geodesic passing through the point $a$ orthogonally
to the diameter $(-1,1)$ and oriented such that $-1\in \partial
H_+$.

  If $E=E_0\cup E_1$ is a compact subset in $\mathbb{D}$ such that $E_0$ is a compact subset
  of $H_+\cup \gamma$ and $E_1\subset [a,1)$ is a compact set such
that $\ell_{\mathbb{D}}(E_1)=\uptau$, then %
\begin{equation} \label{Equation 3.9}
{\rm cap}(E_0\cup[a,\rho(a,\uptau)])\le {\rm cap}(E). %
\end{equation} %

If $E_0$ has positive logarithmic capacity, then equality occurs
in (\ref{Equation 3.9})  if and only if $E_1$ coincides with the
interval $[a,\rho(a,\uptau)]$ up to a set of zero logarithmic
capacity. Otherwise,   equality occurs in (\ref{Equation 3.9}) if
and only if $E_1$ coincides with some interval
$[b,c]\subset [a,1)$ such that $\ell_\mathbb{D}([b,c])=\uptau$ up to a set of zero logarithmic capacity.%
\end{corollary} %

\noindent %
\emph{Proof.} Since the considered characteristics of compact sets
are invariant under M\"{o}bius transformations, we may assume once
more that $a=0$. Notice that the polarization transformations used
in the proof of Lemma~\ref{Lemma 3.4} do not change the portion
$E_0$ of the set $E$, which  lies  in the halfspace $H_+=\{z\in
\mathbb{D}:\, \Re z\le 0\}$. Therefore, our arguments used in the
proof of Lemma~\ref{Lemma 3.4} prove this corollary as well.
 \hfill $\Box$

\medskip

As concerns upper bounds for the conformal capacity of compact
sets $E\subset (-1,1)$ having fixed hyperbolic length, it is
expected that there are no non-trivial upper bounds in this case.
Below we present two examples which confirm these expectations.

\begin{example}\label{Example 3.1}
It was shown by M.~Tsuji \cite{Tsuji-1} that the standard Cantor
set $K$ has positive logarithmic capacity; more precisely, ${\rm
log.cap}(K)\ge 1/9$, see \cite[p. 143]{R}. Hence, by
Proposition~\ref{Proposition 2.2.1}, the conformal capacity
$\kappa={\rm cap}(K_{1/3})$ of the scaled Cantor set
$K_{1/3}=\{z:\,3z\in K\}$ is positive.

For $n\in \mathbb{N}$, let $K_{1/3}^n=\varphi_n(K_{1/3})$ denote
the image of the scaled Cantor set $K_{1/3}$ under the M\"{o}bius
mapping $\varphi_n(z)=(z+r_n)/(1+r_nz)$, where
$r_n=1-\frac{1}{n!}$. Since the hyperbolic length and conformal
capacity are invariant under M\"{o}bius automorphisms of
$\mathbb{D}$, we have
$\ell_\mathbb{D}(K_{1/3}^n)=\ell_\mathbb{D}(K_{1/3})=0$ and ${\rm
cap}(K_{1/3}^n)={\rm cap}(K_{1/3})=\kappa$ for all $n\in
\mathbb{N}$.
A simple calculation shows that %
$$ %
p_{\mathbb{D}}(r_{n+1},r_n)=\frac{n}{n+2-\frac{1}{n!}}, \quad n\in
\mathbb{N}.
$$ %
This implies that the sequence of pseudo-hyperbolic distances
$p_\mathbb{D}(r_{n+1},r_n)$ strictly increases and
$p_\mathbb{D}(r_{n+1},r_n)\to 1$ as $n\to \infty$. Therefore, the
sequence of hyperbolic distances $d_\mathbb{D}(r_{n+1},r_n)$ also
strictly increases and $d_\mathbb{D}(r_{n+1},r_n)\to \infty$ as
$n\to \infty$. This implies that $K_{1/3}^k$ and $K_{1/3}^l$ are
disjoint if $k\not=l$ and %
\begin{equation}   \label{Equation 3.6.1} %
d_\mathbb{D}(K_{1/3}^{n+1},K_{1/3}^n)\to \infty \quad {\mbox{ as $n\to \infty$.}} %
\end{equation} %

Given any $C>0$, we fix $j\in \mathbb{N}$ such that $j\kappa>C$.
Then, for any $m\in \mathbb{N}$, we consider the compact sets
$K_{1/3}^{m,j}=\cup_{s=m}^{m+j-1}K_{1/3}^s$. Using (\ref{Equation
3.6.1}) and arguing as in the proof of the limit relation
(\ref{Equation 2.20}) of Proposition~\ref{Proposition 2.3.2}, we
conclude that ${\rm cap}(K_{1/3}^{m,j})\to j\kappa$ as $m\to
\infty$.  The latter shows that for every constant $C>0$ there are
compact sets $E\subset (-1,1)$ such
that $\ell_\mathbb{D}(E)=0$ and ${\rm cap}(E)\ge C$. %
\end{example} %


\medskip

In our previous example, the hyperbolic diameters of sets
$K_{1/3}^{m,j}$ tend to $\infty$ as $m\to \infty$. For compact
sets with hyperbolic diameters bounded by some constant, say for
compact sets on the interval $[a,b]\subset (-1,1)$ such that
$\ell_\mathbb{D}(E)<\ell_\mathbb{D}([a,b])$, we have the strict
inequality ${\rm cap}(E)<{\rm cap}([a,b])$, which follows from the
fact that $[a,b]\setminus E$ contains a non-empty open interval
that is a set of positive logarithmic capacity. In our next
example, we show that
for every $\varepsilon>0$ and $a,b$ and $\uptau$ are such that
$-1<a<b<1$, $0<\uptau<\ell_{\mathbb{D}}([a,b])$, there is a
compact set $E\subset [a,b]$ such that $\ell_{\mathbb
D}(E)=\uptau$ and ${\rm cap}(E)>{\rm cap}([a,b])-\varepsilon$.

\begin{example}\label{Example 3.2} %
First, we consider a condenser $(A(\rho^{-1},\rho),K_n(l))$ with
the domain $A(\rho^{-1},\rho)$, where
$A(\rho_1,\rho_2)=\{z:\,\rho_1<|z|<\rho_2\}$, $0<\rho_1<\rho_2$,
and a compact set $K_n(l)=\cup_{j=1}^n K_{n,j}(l)$, where
$K_{n,j}(l)=\{e^{i\theta}:\,|\theta-2\pi(j-1)/n|\le l/n\}$,
$j=1,\ldots,n$, $0<l<\pi$.

Let $\Gamma=\Gamma(\rho,n,l)$ denote the family of curves in
$A(\rho^{-1},\rho)$ joining the boundary circles of
$A(\rho^{-1},\rho)$ with the set $K_n(l)$ and let
$\Gamma_0=\Gamma_0(\rho,n,l)$ denote the family of curves in the
annular sector $S_n(\rho)=\{z:\,1<|z|<\rho,\ |\arg z|<\pi/n\}$
joining the arc $\{z=\rho e^{i\theta}:\,|\theta|\le \pi/n\}$ with
the set $K_{n,1}(l)$. It follows from Ziemer's relation between
the capacity of a condenser and the modulus of an appropriate
family of curves (see \cite[Theorem 3.8]{Zim}) and from symmetry
properties of the
modulus of family of curves (see \cite[Theorem 4]{Ahlfors2006} for an equivalent form of this symmetry property given in terms of the extremal length)that %
\begin{equation} \label{Equation 3.12}
{\rm cap}((A(\rho^{-1},\rho),K_n(l)))=\M(\Gamma)=2n \M(\Gamma_0).
\end{equation} %

To find $\M(\Gamma_0)$, we consider a function
$\varphi(z)=\varphi_2\circ \varphi_1(z)$ with %
\begin{equation} \label{Equation 3.13}%
\varphi_1(z)=i(n \mathcal{K}(k)/\pi)\log z \quad {\mbox{and}}
\quad \varphi_2(\zeta)=\sn(\zeta,k),
\end{equation} %
where the parameter $k$, $0<k<1$, of the elliptic sine function is
defined by the equation %
\begin{equation} \label{Equation 3.14} %
\frac{\mathcal{K}'(k)}{\mathcal{K}(k)}=\frac{2n\log\rho}{\pi}.
\end{equation} %
It is a well-known property of elliptic integrals \cite[Theorem
5.13(1)]{AVV} that if $\mathcal{K}'(k)/\mathcal{K}(k)\to \infty$,
then $k\to 0$ and the
following expansion holds: %
\begin{equation} \label{Equation 3.15} %
\frac{\mathcal{K}'(k)}{\mathcal{K}(k)}=\frac{2}{\pi}\log\frac{4}{k}+o(1).
\end{equation} %
From (\ref{Equation 3.14}) and (\ref{Equation 3.15}), we obtain
that %
\begin{equation} \label{Equation 3.16}
\log\frac{1}{k}=n\log \rho-\log 4+o(1),
\end{equation} %
where $o(1)\to 0$ when $n\to \infty$.

\smallskip

The function $w=\varphi(z)$ maps $S_n(\rho)$ conformally onto the
semidisk $\mathbb{D}_R^+=\{w\in \mathbb{D}_R:\,\Im w>0\}$ with
$R=1/\sqrt{k}$. Furthermore, this function maps the arc
$K_{n,1}(l)$ onto an interval $[-c_n(l),c_n(l)]$ with
$0<c_n(l)<1$. To find $c_n(l)$, we note that $\mathcal{K}(k)\to
\pi/2$ as $k\to 0$ and that $\sn(\zeta,k)$ converges to
$\sin\zeta$ uniformly on compact subsets of $\mathbb{C}$ as $k\to
0$. Using these relations and equations (\ref{Equation 3.13}), we
find that $c_n(l)=\sin(l/2)+o(1)$ and, therefore, we have the
following asymptotic formula for the logarithmic capacity of the
interval $[-c_n(l),c_n(l)]$:
\begin{equation} \label{Equation 3.17} %
{\rm {log.cap}}([-c_n(l),c_n(l)])=
(1/2)\sin(l/2)+o(1),
\end{equation} %
where $o(1)\to 0$ as $n\to\infty$.

Let $\Gamma_1=\Gamma_1(\rho,n,l)$ denote the family of curves in
$\mathbb{D}_R$ joining the circle $\mathbb{T}_R$ with the interval
$[-c_n(l),c_n(l)]$. Conformal invariance and symmetry properties
of the modulus of family of curves  imply that %
\begin{equation} \label{Equation 3.18} %
\M(\Gamma_0)=\frac{1}{2}\M(\Gamma_1).
\end{equation} %
We have the following limit relation between the modulus of
$\Gamma_1$ and the logarithmic capacity of
$[-c_n(l),c_n(l)]$: %
\begin{equation} \label{Equation 3.19} %
-\frac{1}{2\pi}\log({\rm
{log.cap}}([-c_n(l),c_n(l)]))=(\M(\Gamma_1))^{-1}-\frac{1}{2\pi}\log
R+o(1),
\end{equation} %
where $o(1)\to 0$ as $R\to\infty$.

Using relations (\ref{Equation 3.17}) and (\ref{Equation 3.19})
with $R=1/\sqrt{k}$ and with $k$ as in the equation
(\ref{Equation 3.16}), we find that %
\begin{equation} \label{Equation 3.20} %
\M(\Gamma_1)=\left(-\frac{1}{2\pi}\log\frac{\sin(l/2)}{2}+\frac{1}{4\pi}(n\log
\rho-\log 4)+o(1)\right)^{-1}.
\end{equation} %

Finally, combining equations (\ref{Equation 3.12}), (\ref{Equation
3.18}) and (\ref{Equation 3.20}), we obtain the following
asymptotic formula for the capacity of the condenser
$(A(\rho^{-1},\rho),K_n(l))$:
\begin{equation} \label{Equation 3.21} %
{\rm cap}((A(\rho^{-1},\rho),K_n(l)))=\frac{4\pi}{\log\rho}+o(1),
\end{equation} %
where $o(1)\to 0$ when $\rho$ and $l$ are fixed and $n\to\infty$.

Let $s=s(\rho)$, $0<s<1$, be such that %
\begin{equation}  \label{Equation 3.22} %
{\rm cap}([-s,s])={\rm cap}((\mathbb{D}_\rho,\overline{\mathbb{D}}))=\frac{2\pi}{\log\rho}. %
\end{equation} %
Then there is a unique function $\psi(z)$ mapping $A(1,\rho)$
conformally onto $\mathbb{D}\setminus [-s,s]$ such that
$\psi(\rho)=1$. This function can be extended to a  function
continuous on $\overline{A(1,\rho)}$ and such that
$\psi(\overline{z}) =\overline{\psi(z)}$ for all $z\in
\overline{A(1,\rho)}$. Thus, $\psi(z)$ maps $K_n(l)$ onto a
compact set $\psi(K_n(l))\subset [-s,s]$.

The same conformal invariance and symmetry properties, which we
used earlier in this example, together with equations (\ref{Equation 3.21}) and (\ref{Equation 3.22}) imply that %
$$ 
{\rm cap}(\psi(K_n(l)))=\frac{1}{2}{\rm cap}((A(\rho^{-1},\rho),K_n(l)))=\frac{2\pi}{\log\rho}+o(1). %
$$ 

Notice that there is a constant $C>0$ such that
$|\psi'(e^{i\theta})|\le C$ for all $\theta\in \mathbb{R}$. This
implies that for every $\uptau>0$ there is $l_0$, $0<l_0<2\pi$,
such that $\ell_\mathbb{D}(\psi(K_n(l)))\le \uptau$ for every $l$,
$0<l\le l_0$, and all $n\ge 2$.

We fix $l$, $0<l\le l_0$, and $n\ge 2$ and consider a compact set
$E\subset [-s,s]$ such that $\psi(K_n(l))\subset E$ and
$\ell_\mathbb{D}(E)=\uptau$, Then, if $n$ is large enough, we have %
$$ 
{\rm cap}([-s,s])>{\rm cap}(E)\ge {\rm cap}(\psi(K_n(l)))=\frac{2\pi}{\log\rho}+o(1). %
$$ 
The latter equation shows that for every $\uptau$,
$0<\uptau<\ell_\mathbb{D}([-s,s])$, and every $\varepsilon>0$,
there is a compact set $E\subset [-s,s]$ such that
$\ell_\mathbb{D}(E)=\uptau$ and ${\rm cap}(E)>{\rm
cap}([-s,s])-\varepsilon$. Since the hyperbolic length and
conformal capacity are invariant under M\"{o}bius automorphisms of
$\mathbb{D}$, compact sets with similar properties exist for
every interval $[a,b]$, $-1<a<b<1$. %
\end{example}

\medskip

\begin{rem}
Our construction of a compact set in Example~\ref{Example 3.2} is
similar to the construction used in \cite{Solynin1996b} to provide
a counterexample for P.M.~Tamrazov's conjecture on the capacity of
a condenser with plates of prescribed transfinite diameters. In
turn, a counterexample used in \cite{Solynin1996b} is based on the
following result, that is example 5) in \cite[Ch. II,
\S4]{Lebedev}:

Let $K_n^+(l)\subset [-1,1]$ denote the orthogonal projection of
the set $K_n(l)$ introduced earlier onto the real axis. Then %
$$ %
{\rm log.cap}(K_n^+(l))=(1/2)\left(\sin(l/2)\right)^{2/n}. %
$$ %
Taking the limit in this equation as $n\to \infty$, we conclude
that for every $0<s<2$ and every $\varepsilon>0$, there exists a
compact set $E\subset [-1,1]$ with
Euclidean length $s$ such that %
$$ %
{\rm log.cap}(E)>{\rm log.cap}([-1,1])
-\varepsilon=1/2-\varepsilon. %
$$ %
In particular, this answers a question raised in Problem 2 in
\cite{Betsakos2005} by showing that the supremum of the
transfinite diameters of compact sets $E\subset [-1,1]$ with
Euclidean length $s$, $0<s<2$, is equal to $1/2$.
\end{rem}

\medskip

As Examples~\ref{Example 3.1} and \ref{Example 3.2} show, there
are no upper bounds for the capacity expressed in terms of the
hyperbolic length of a compact set $E$, in general. In a
particular case, when $E$ is connected, a non-trivial upper bound
exists and is given in the following lemma. %

\begin{lemma} \label{Lemma 3.3}
Let $L\subset \mathbb{D}$ be a Jordan arc having the hyperbolic length $\uptau>0$. Then %
\begin{equation}  \label{Equation 3.25} %
{\rm cap}(L)\le {\rm cap}([0,r(\uptau)])=4\frac{\mathcal{K}((e^\uptau-1)/(e^\uptau+1))}{\mathcal{K}'((e^\uptau-1)/(e^\uptau+1))}.%
\end{equation}  %
\end{lemma} %
\smallskip

\noindent %
\emph{Proof.} The proof repeats the well-known proof for the
logarithmic capacity, see, for example, \cite[Theorem 5.3.2]{R}.
We consider a parametrization $T:\,[0,r(\uptau)] \to L$ of $L$ by
the hyperbolic arc-length. Then $T$ is contractive in the
hyperbolic metric. Therefore, (\ref{Equation 3.25}) follows from
Proposition~\ref{Proposition 2.19}. \hfill $\Box$

\medskip

 The polarization technique used in
the proof of Lemma~\ref{Lemma 3.1} can be applied in a more
general situation as we demonstrate in our next theorem.


\begin{theorem} \label{Theorem 3.26}
    Consider $1\le n\le 4$ distinct radial intervals $I_k=I(\alpha_k)$, $k=1,\ldots,n$, of
    the unit disk  $\mathbb D$. Suppose that $E_k\subset I_k$,
    $k=1,\ldots, n$,  is a compact set on $I_k$ such that
    $\ell_{\mathbb{D}}(E_k)=l_k$ and that $E_k^*\subset I_k$ is a
    hyperbolic interval having one end point at $z=0$ such that
    $\ell_{\mathbb{D}}(E_k^*)=l_k$.

    If each of the angles formed by the radial intervals $I_k$ and
    $I_j$, $k\not= j$, is greater than or equal to $\pi/2$, then
    \begin{equation}\label{Equation 3.27}
    {\rm cap}\left(\cup_{k=1}^n E_k\right)\ge {\rm cap}\left(\cup_{k=1}^n E_k^*\right).
    \end{equation} %
    Equality occurs in (\ref{Equation 3.27}) if and only if for each $k$,
    $E_k$ coincides with $E_k^*$ up to a set of zero logarithmic capacity.
\end{theorem}

\smallskip

\noindent %
\emph{Proof.} The proof is the same for all $n$. Thus, we assume
that $n=4$. Rotating, if necessary, we may assume that
$I_1=[0,1)$. The diameter $(-i,i)$ is a hyperbolic geodesic, which
we orient such that $1/2\in H_+$. If $n=4$, then the angles
between the neighboring intervals are equal to $\pi/2$ and
therefore the set $K=\cup_{k=1}^4 E_k$ can be represented as the
union $K=E_0\cup E_1$ with $E_0=\cup_{k=2}^4 E_k\subset
\overline{H_+}$. This shows that the sets $E_0$, $E_1$ and
$K_1=E_1^*\cup E_2\cup E_3\cup E_4$ satisfy the assumptions of
Corollary~\ref{Corollary 3.1}.
Therefore, by this corollary,%
$$ %
{\rm cap}(K_1)={\rm cap}(E_0\cup E_1^*)\le {\rm cap}(E_0\cup
E_1)={\rm cap}(K)
$$ %
with the sign of equality if and only if $E_1$ coincides with
$E_1^*$ up to a set of zero logarithmic capacity.

The same argument can be applied successively to the sets $K_1$,
$K_2=E_1^*\cup E_2^*\cup E_3\cup E_4$, $K_3=E_1^*\cup E_2^*\cup
E_3^*\cup E_4$ and $K_4=E_1^*\cup E_2^*\cup
E_3^*\cup E_4^*$ to obtain the inequalities %
$$ %
{\rm cap}\left(\cup_{k=1}^4 E_k^*\right)={\rm cap}(K_4)\le {\rm
cap}(K_3)\le {\rm cap}(K_2)\le {\rm cap}(K_1)\le{\rm cap}(K)={\rm
cap}\left(\cup_{k=1}^4 E_k\right).
$$ %
Moreover, equality occurs in any one of these inequalities if and
only if the corresponding sets $E_k$ and $E_k^*$ coincide up to a
set of zero logarithmic capacity. Thus, the theorem is proved.
 \hfill $\Box$

\medskip

\begin{rem} %
 We want to stress once more that the
inequality (\ref{Equation 3.27}) also follows from the contraction
principle of Proposition~\ref{Proposition 2.19}. Indeed, under the
assumptions that all angles between radial intervals are $\ge
\pi/2$, there is always a
contraction $\varphi:\cup_{k=1}^n E_k\to \cup_{k=1}^n E_k^*$. %
\end{rem}

\medskip



If the angle between some radial intervals $I_k$ and $I_j$ is
smaller than $\pi/2$, then both our proofs, with polarization or
with the contraction principle, fail even in the simplest case of
two intervals $I_1$ and $I_2$ and when each of the sets
$E_1\subset I_1$ and $E_2\subset I_2$ is a hyperbolic interval.
However, the graphs of the results of numerical experiments
performed by Dr. Mohamed Nasser, which are presented in Figure~3,
suggest that the monotonicity property of the conformal capacity
of two intervals remains in place for all angles. Therefore, we
suggest the following.

\begin{problem} %
Given fixed $0<r<1$, $0<s<\infty$, and $0<\alpha<\pi/2$ and
varying $0\le t<1$, let $E(t)=[0,r]\cup \{\uptau e^{i\alpha}: t\le
\uptau \le d(s,t)\}$ with $d(s,t)$ such that
$\ell_{\mathbb{D}}([t,d(s,t)])=s$. Prove (or disprove) that ${\rm
cap}(E(t))$ strictly increases on the interval $0\le t<1$.
\end{problem}

\begin{figure} %
\hspace{-0.5cm} %
\begin{minipage}{1.0\textwidth}
$$\includegraphics[scale=0.55,angle=0]{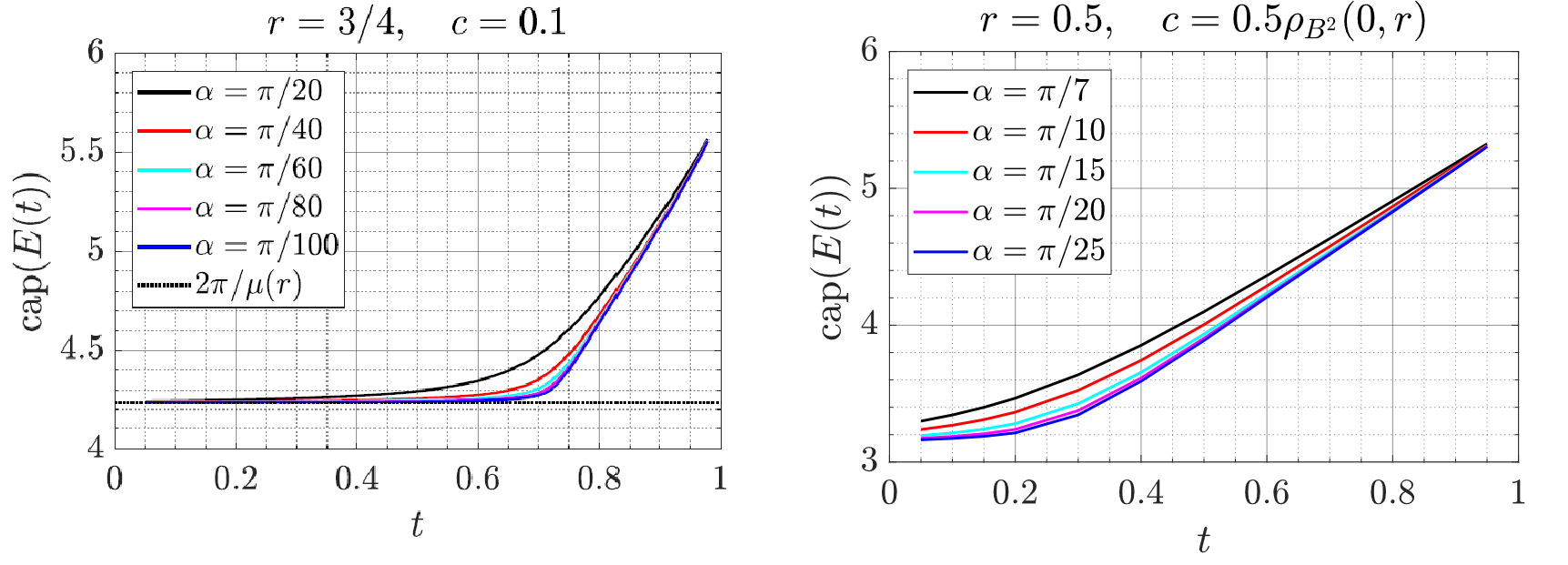}
$$
\end{minipage}
\caption{Hyperbolic capacity of two intervals.} %
\end{figure}

\medskip

As one can see from our proof of Theorem~\ref{Theorem 3.26}, the
restriction $n\le 4$ on the number of radial intervals and
restriction on angles between them is needed because otherwise
polarization with respect to hyperbolic geodesics may destroy the
radial structure of compact sets under consideration. Also, if at
least one angle between radial intervals is $<\pi/2$, then the
contraction principle of Proposition~\ref{Proposition 2.19} can
not be applied, in general. Still, under some additional
assumptions on the hyperbolic lengths and angles, we have the
following more general version.

    \begin{theorem} \label{Theorem 3.30}%
        Let $E_k$, $k=1,\ldots,n$, be compact sets on the radial intervals
        $I_k=[0,e^{i\beta_k})$,    $0=\beta_1<\beta_2<\ldots<\beta_n<\beta_{n+1}=2\pi$, such that %
       $$ %
        \ell_{\mathbb{D}}(E_k)\ge 2 \log\cot\frac{\alpha}{2}, \quad k=1,\ldots,n, %
        $$ %
        where $\alpha$ stands for the minimal angle between the intervals
        $I_k$. Then
        \begin{equation}\label{Equation 3.31}
        {\rm cap}\left(\cup_{k=1}^n E_k\right)\ge {\rm cap}
        \left(\cup_{k=1}^n E_k^*\right),
        \end{equation} %
        where  $E_k^*\subset
        I_k$ is a hyperbolic interval having one end point at $z=0$ such
        that $\ell_{\mathbb{D}}(E_k^*)=\ell_{\mathbb{D}}(E_k),    k=1,\ldots,n$.

        Equality occurs here if and only if for each $k$, $E_k$ coincides
        with $E_k^*$ up to a set of zero logarithmic capacity.
    \end{theorem}

    \smallskip

    \noindent %
    \emph{Proof.} (a) Let $\uptau_\alpha=\log\cot(\alpha/2)$,
    $r_\alpha=r(\uptau_\alpha)=\sec\alpha-\tan\alpha$ and let
    $\gamma_\alpha$ be the hyperbolic geodesic orthogonal to
    the interval $(-1,1)$ at $z=r_\alpha$ and oriented such that
    $0$ belongs to the halfplane $H_+=H_+(\gamma_\alpha)$. An easy
    calculation shows that $\gamma_\alpha$ has its endpoints at
    the points $e^{\pm i\alpha}$. Consider sets $E_1^+=E_1\cap
    \overline{H_+}$, $E_1^-=E_1\setminus H_+$ and
    $E_0=E_1^+\cup\left(\cup_{k=2}^n E_k\right)$. Let $\widetilde{E}_1^-$
    denote the closed hyperbolic interval with the initial point at $z=r(\alpha)$ such that
    $\widetilde{E}_1^-\subset [r_\alpha,1)$ and $\ell_\mathbb{D}(\widetilde{E}_1^-)=\ell_\mathbb{D}(E_1^-)$
    and let $\widetilde{E}_1=E_1^+\cup \widetilde{E}_1^-$. Since the
    minimal angle between the intervals $I_k$ is $\alpha$ and
    $\gamma_\alpha$ has its endpoints at $e^{\pm i\alpha}$, it
    follows that $E_0\subset \overline{H_+}$. Therefore, we can
    apply Corollary~\ref{Corollary 3.1} to obtain the following: %
    \begin{equation}  \label{Equation 3.32} %
{\rm cap}(E_0\cup\widetilde{E}_1^-)={\rm
cap}(\widetilde{E}_1\cup(\cup_{k=2}^n E_k))\le {\rm cap}(E)
    \end{equation} %
with the sign of equality if and only if $E_1^-$ coincides with
$\widetilde{E}_1^-$ up to a set of zero logarithmic capacity.

If $E_1^+=[0,r_\alpha]$, then $E_1^*=E_1^+\cup \widetilde{E}_1^-$
and
the inequality in (\ref{Equation 3.32}) is equivalent to the inequality %
\begin{equation} \label{Equation 3.33} %
{\rm cap}\left(E_1^*\cup\left(\cup_{k=2}^n E_k\right)\right)\le
{\rm cap}(E)
\end{equation} %
with the sign of equality if and only if $E_1^*$ coincides with
$E_1$ up to a set of zero logarithmic capacity.

\smallskip

(b) If $E_1^+\not=[0,r_\alpha]$, then $[0,r_\alpha]\setminus
E_1^+$ contains an open interval. In this case
$\ell_\mathbb{D}(\widetilde{E}_1^-)>\ell_\mathbb{D}([0,r_\alpha])$.
Let $b_1$, $r_\alpha<b_1<1$,  denote the end point of the interval
$\widetilde{E}_1^-$ and let $c_1$ denote the midpoint of the
hyperbolic interval $[0,b_1]_h$. Notice that under our assumptions
$r_\alpha<c_1$. Let $\gamma_1$ be the hyperbolic geodesic
orthogonal to the diameter $(-1,1)$ at $z=c_1$ and oriented such
that $0\in H_+(\gamma_1)$. Let
$\widehat{E}_1=\mathcal{P}_{\gamma_1}(\widetilde{E}_1)$ and
$\widehat{E}=\mathcal{P}_{\gamma_1}(\widetilde{E}_1\cup(\cup(_{k=2}^n(E_k)))$
denote polarizations of the corresponding sets with respect to
$\gamma_1$. Since $r_\alpha<c_1$ and therefore $\cup_{k=2}^n
E_k\subset H_+(\gamma_1)$, we have $\widehat{E}_1\supset
[0,r_\alpha]$ and $\widehat{E}=\widehat{E}_1\cup(\cup_{k=1}^n
E_k)$. Applying Proposition~\ref{Proposition 2.5} and using
(\ref{Equation 3.32}), we conclude
that %
\begin{equation} \label{Equation 3.34} %
{\rm cap}(\widehat{E}_1\cup(\cup_{k=2}^n E_k))<{\rm
cap}(\widetilde{E}_1\cup(\cup_{k=2}^n E_k))\le {\rm cap}(E)
\end{equation} %
with the sign of strict inequality in the first inequality because
$\widehat{E}_1$ differs from $\widetilde{E}_1$ and from its
reflection with respect to $\gamma_1$ by an open interval and
therefore by a set of positive logarithmic capacity.

Since $\ell_\mathbb{D}(\widehat{E}_1)=\ell_\mathbb{D}(E_1)$ and
$\widehat{E}_1\supset [0,r_\alpha]$, using (\ref{Equation 3.34})
and applying the same arguments as in part (a) of this proof to
the set $\widehat{E}=\widehat{E}_1\cup(\cup_{k=2}^n)$, we conclude
that in the case $E_1^+\not=[0,r_\alpha]$ the inequality
(\ref{Equation 3.33}) remains true with the same statement on the
equality cases.

\smallskip

(c) Now, when (\ref{Equation 3.33}) is proved in all cases, we can
apply the iterative procedure as in the proof of
Theorem~\ref{Theorem 3.26} to conclude that the inequality
(\ref{Equation 3.31}) holds with  the sign of equality if and only
if, for each $k$, $E_k$ coincides  with $E_k^*$ up to a set of
zero
logarithmic capacity.  \hfill  $\Box$ %

\medskip

For compact sets $E_1$ and $E_2$,  lying on two orthogonal
diameters of $\mathbb{D}$, we have the following result.

\begin{theorem} \label{Theorem 3.35}%
Let $E_1\subset (-1,1)$, $E_2\subset (-i,i)$ be compact sets and
 let $r_k$, $0<r_k<1$, and $l_k>0$ be such that %
$$ %
\ell_\mathbb{D}(E_k)=\ell_\mathbb{D}([-r_k,r_k])=2l_k, \quad
k=1,2.
$$ %
Then %
\begin{equation} \label{Equation 3.36}%
{\rm cap}(E_1\cup E_2)\ge {\rm cap}([-r_1,r_1]\cup
[-ir_2,ir_2])\ge {\rm cap}([-r_0,r_0]\cup [-ir_0,ir_0]),
\end{equation}  %
where $0<r_0<1$ is such that
$\ell_\mathbb{D}([-r_0,r_0])=l_1+l_2$.


 Equality occurs in the first inequality if and only if $E_1$ coincides
 with $[-r_1,r_1]$ and $E_2$ coincides
 with $[-ir_2,ir_2]$ up to a set of zero logarithmic capacity.
Equality occurs in the second inequality if and only if
$l_1=l_2$.
\end{theorem}

\smallskip

\noindent %
\emph{Proof.} Let $0\le r_k^{\pm}<1$, $k=1,2$, be such that the
following holds:
$$ %
\ell_\mathbb{D}([0,r_1^+])=\ell_\mathbb{D}(E_1\cap [0,1)), \quad
\ell_\mathbb{D}([0,-r_1^-])=\ell_\mathbb{D}(E_1\cap [0,-1)), %
$$
$$
\ell_\mathbb{D}([0,r_2^+])=\ell_\mathbb{D}(E_2\cap [0,i)), \quad
\ell_\mathbb{D}([0,-r_2^-])=\ell_\mathbb{D}(E_2\cap [0,-i)).
$$ %
Then, by Theorem~\ref{Theorem 3.26},
\begin{equation} \label{Equation 3.37}%
{\rm cap}(E_1\cup E_2)\ge {\rm cap}([0,r_1^+]\cup [0,-r_1^-]\cup
[0,ir_2^+)\cup [0,-ir_2^-])
\end{equation}  %
with the sign of equality if and only if the sets in the left and
right sides of this inequality coincide up to a set of zero
logarithmic capacity.

Suppose that $r_1^+\not= r_1^-$, say $r_1^+>r_1^-$.  Then
$r_1^+>r_1>r_1^-$. Let $\mathcal{P}_\gamma$ denote polarization
with respect to the geodesic $\gamma$ that is orthogonal to the
hyperbolic interval $[-r_1,r_1^+]$ at its midpoint $c$, $0<c<r_1$.
We assume here that $\gamma$ is oriented such that $0\in
H_\gamma^+$. Under our assumptions,
$\mathcal{P}_\gamma([-r_1^-,r_1^+]\cup[-ir_2^-,ir_2^+])=[-r_1,r_1]\cup[-ir_2^-,ir_2^+]$.
Since the set $[-r_1^-,r_1^+]\setminus [-r_1,r_1]$ has positive
logarithmic capacity, it follows from
Proposition~\ref{Proposition 2.5} that %
\begin{equation} \label{Equation 3.38} %
{\rm cap}([-r_1^-,r_1^+]\cup [-ir_2^-,ir_2^+])>{\rm
cap}([-r_1,r_1]\cup [-ir_2^-,ir_2^+]).
\end{equation}  %

The same polarization argument can be applied to show that, if
$r_2^-\not= r_2^+$, then %
\begin{equation} \label{Equation 3.39}%
{\rm cap}([-r_1,r_1]\cup [-ir_2^-,ir_2^+])>{\rm
cap}([-r_1,r_1]\cup [-ir_2,ir_2]).
\end{equation}  %
Combining inequalities (\ref{Equation 3.37})--(\ref{Equation
3.39}), we obtain the first inequality in (\ref{Equation 3.36})
with the sign of equality if and only if $E_1$ coincides
 with $[-r_1,r_1]$ and $E_2$ coincides
 with $[-ir_2,ir_2]$ up to a set of zero logarithmic capacity.

\smallskip

To prove the second inequality in (\ref{Equation 3.36}), we use the conformal mapping 
$$ %
g(z)=\sqrt{(z^2+r_2^2)/(1+r_2^2z^2)}
$$ %
from the doubly connected domain $\mathbb{D}\setminus
([-r_1,r_1]\cup [-ir_2,ir_2])$ onto 
$\mathbb{D}\setminus ([-r,r])$  with %
$$ %
 r=\sqrt{(r_1^2+r_2^2)/(1+r_1^2r_2^2)}.
$$ %
We note that conformal mappings preserve capacity of condensers
and that the function ${\rm cap}([-r,r])={\rm
cap}([-r_1,r_1]\cup[-ir_2,ir_2])$ is strictly increasing on the
interval $0\le r<1$. Furthermore, it follows from formulas
(\ref{Equation 1.1}) that the sum $l_1+l_2$ of the hyperbolic
lengths defined in the theorem is constant if and only if the
following product is constant:
\begin{equation} \label{Equation 3.40}
\frac{1+r_1}{1-r_1}\cdot \frac{1+r_2}{1-r_2}=C,
\end{equation} 
where $C$ is constant.

 Our goal now is to minimise the
function $F=F(r_1,r_2)$, defined by%
$$
F=(r_1^2+r_2^2)/(1+r_1^2r_2^2),
$$ %
under the constraint (\ref{Equation 3.40}). %

Introducing new variables %
$ u=(1+r_1)/(1-r_1)$, $v=(1+r_2)/(1-r_2)$ and $w=u^2+v^2$, we can
express $F$ in terms of these variables as follows:
$$ %
F=\frac{(u^2+v^2)+(1-4C+C^2)}{(u^2+v^2)+(1+4C+C^2)}=\frac{w+(1-4C+C^2)}{w+(1+4C+C^2)},
$$ %
which we have to minimize $F$ under the constraint $uv=C$. Differentiating, we find that %
$$ %
\frac{d}{dw}F=\frac{4C}{(w+(1+4C+C^2))^2}>0.
$$ %
Therefore, $F$ takes its minimal value when $w=u^2+y^2$ is as
small as possible. By the classical arithmetic-geometric mean
inequality $u^2+v^2> 2uv=2C$, unless $u=v$. Therefore, the minimal
value of $F$ under the constraint $uv=C$ occurs when $u=v$. The
latter implies that ${\rm cap}([-r_1,r_1]\cup [-ir_2,ir_2])\ge
{\rm cap}([-r_0,r_0]\cup[-ir_0,ir_0])$ with the sign of equality
if and only if $r_1=r_2$. This proves the second inequality in
(\ref{Equation 3.36}).
 \hfill $\Box$

\medskip

\begin{rem} %
 The inequality obtained in Theorem~\ref{Theorem 3.35}
 is stronger, in general, than the inequality obtained by the classical
Steiner symmetrization, which will be discussed in Section~5. To
give an example, we consider two sets:  $E_1=[0,1/2]\cup [0,i/2]$
and $E_2=[0,1/2]\cup [0,i/4]$. Two Steiner symmetrizations, chosen
appropriately, transform these sets to the sets
$E_1^*=[-1/4,1/4]\cup [-i/4,i/4]$ and $E_2^*=[-1/4,1/4]\cup
[-i/8,i/8]$, respectively. The first inequality in (\ref{Equation
3.36}) compares
the conformal capacities of $E_1$ and $E_2$ with the conformal capacities of the sets %
$$  %
E_1^o=[-(2-\sqrt{3}),2-\sqrt{3}]\cup
[-i(2-\sqrt{3}),i(2-\sqrt{3})], %
$$ %
$$ %
E_2^o=[-(2-\sqrt{3}),2-\sqrt{3}]\cup
[-i(4-\sqrt{15}),i(4-\sqrt{15})]. %
$$
Numerical computation gives the following approximation and bounds
for ${\rm cap}(E_k)$, $k=1,2$, %
$$ %
{\rm cap}(E_1^*)\approx 3.62589<{\rm cap}(E_1^o)\approx
3.77702<{\rm cap}(E_1)\approx 4.28254,
$$ %
$$ %
{\rm cap}(E_2^*)\approx 3.19333< {\rm cap}(E_2^o)\approx 3.29244
<{\rm cap}(E_2)\approx 3.60548.
$$ %
\end{rem} %

\medskip

In our previous lemmas and theorems of this section, the extremal
configurations were hedgehogs with spikes issuing from the central
point. Similar results for compact sets with spikes emanating
from a certain compact central body $E_0$ sitting in the disk
$\overline{\mathbb{D}}_r$, $0<r<1$, as it is shown in Figure~4,
also may be useful in applications.

\begin{theorem} \label{Theorem 3.42}%
        Let $E_0$ be a compact set in the disk $\overline{\mathbb{D}}_r$, $0<r<1$.
        Let $E_k$, $k=1,\ldots,n$, be compact sets on the radial intervals
        $I_k=[re^{i\beta_k},e^{i\beta_k}]$,
        $0=\beta_1<\beta_2<\ldots<\beta_n<\beta_{n+1}=2\pi$, such that %
        $$ 
        \ell_{\mathbb{D}}(E_k)\ge 2 \log \left(\cot\frac{\alpha}{2}\frac{1-r}{1+r}\right), \quad
        k=1,\ldots,n,
        $$ 
        where $\alpha$ stands for the minimal angle between the intervals
        $I_k$. Then
        \begin{equation}\label{Equation 3.43}
        {\rm cap}(\cup_{k=0}^n E_k)\ge {\rm cap}(\cup_{k=0}^n
        E_k^*),
        \end{equation} %
        where  $E_0^*=E_0$ and $E_k^*\subset
        I_k$ is a hyperbolic interval having one end point at $z=re^{i\beta_k}$ such
        that $\ell_{\mathbb{D}}(E_k^*)=\ell_{\mathbb{D}}(E_k)$, $k=1,\ldots,n$.

        Equality occurs in (\ref{Equation 3.43}) if and only if for each $k=1,\ldots,n$, $E_k$ coincides
        with $E_k^*$ up to a set of zero logarithmic capacity.

    \end{theorem}

    \begin{figure} 
\begin{minipage}{1.0\textwidth}
$$\includegraphics[scale=0.5,angle=0]{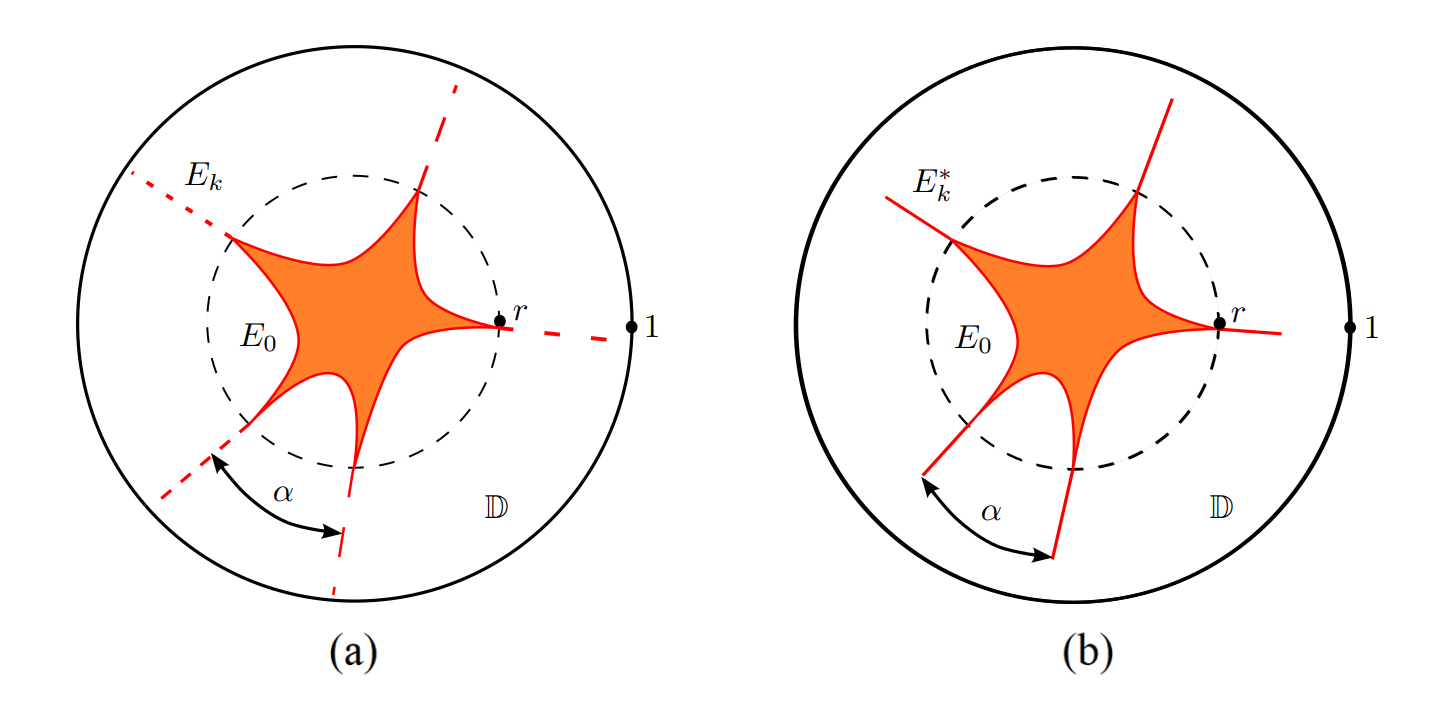}
$$
\end{minipage}
\caption{Hedgehogs with central body $E_0$ and spikes on five intervals.}%
\end{figure}


\noindent %
\emph{Proof.} The proof of this theorem is essentially the same as
the proof of Theorem~\ref{Theorem 3.30}. The only new thing we
need is the following observation. Let $\rho=\rho(r,\alpha)$,
$r<\rho<1$, be
such that %
$$
\ell_\mathbb{D}([r,\rho])=2 \log
\left(\cot\frac{\alpha}{2}\frac{1-r}{1+r}\right) %
$$ %
and let $c$ be the middle point of the hyperbolic interval
$[r,\rho]_h$. Then every polarization $\mathcal{P}_\gamma$
performed with respect to the geodesic $\gamma$ that is orthogonal
to the radial interval $[0,e^{i\beta_k})$ and crosses it at some
point $c_1$ such that $c<|c_1|<1$, and such that $0\in
H_+(\gamma)$, satisfies the following properties: %
$$ %
\mathcal{P}_\gamma(E_k)\subset [re^{i\beta_k},e^{i\beta_k}) \quad
{\mbox{and}} \quad \mathcal{P}_\gamma(E_j)=E_j \quad {\mbox{if
$j\not=k$.}} %
$$ %
Therefore, applying polarizations successively and arguing as in
the proof of Theorems~\ref{Theorem 3.26} and \ref{Theorem 3.30},
we obtain inequality (\ref{Equation 3.43}) together
with the statement on the equality in it. %
\hfill  $\Box$

\medskip

Most of the results presented above in this section can be proved
by two methods, either using the polarization technique or the
contraction principle. Now, we give an example of a result, when
polarization does not work but the contraction principle is easily
applicable.

\begin{theorem} \label{Theorem 3.44}%
Let $E_0$ be a compact subset of $\mathbb{T}=\partial \mathbb{D}$
and let $\rho(\theta)\ge 0$ be an upper semicontinuous function on
$E_0$. For $0<r<1$, let $E(r)$ be a compact set in $\mathbb{D}$
such that the intersection $E(r)\cap [0,e^{i\theta}]$ is empty if
$e^{i\theta}\not\in E_0$ and it is an interval
$[re^{i\theta},se^{i\theta}]$  with $s$, $r\le s<1$,  such that
$\ell_{\mathbb{D}}([r,s])=\rho(\theta)$, if $e^{i\theta}\in E_0$.

Then the conformal capacity ${\rm cap}(E(r))$ is an increasing
function on the interval $0\le r<1$.
\end{theorem}

\smallskip

\noindent %
\emph{Proof.} Let $r_0$, $0<r_0<1$, be fixed and let
$\widetilde{A}(r_0)=\{z:\,r_0\le |z|<1\}$. Consider the function
$\varphi_r:\widetilde{A}(r_0)\to \mathbb{D}$ defined as follows:
if $z=se^{i\theta}\in \widetilde{A}(r_0)$, then $\arg
\varphi_r(z)=\theta$ and
$d_\mathbb{D}(re^{i\theta},\varphi_r(z))=d_\mathbb{D}(r_0e^{i\theta},z)$.
We claim that, if $0\le r\le r_0$, then $\varphi_r$ is a
hyperbolic contraction on $\widetilde{A}(r_0)$.

To prove this claim, we fix two points $z_1,z_2\in
\widetilde{A}(r_0)$ and consider the hyperbolic distance
$d_\mathbb{D}(\varphi_r(z_1),\varphi_r(z_2))$. Our claim will be
proved if we show that this hyperbolic distance or, equivalently,
the pseudo-hyperbolic distance %
\begin{equation} \label{Equation 3.45}%
p_\mathbb{D}(\varphi_r(z_1),\varphi_r(z_2))=\left|\frac{\varphi_r(z_1)-\varphi_r(z_2)}{1-\varphi_r(z_1)\overline{\varphi_r(z_2)}}\right|,
\end{equation} %
is an increasing function of $r$ on the interval $0\le r\le r_0$.

Rotating, if necessary, we may assume that %
$$  %
z_1=a_0e^{i\alpha},\ \ z_2=b_0e^{-i\alpha}, \quad {\mbox{where
$r_0\le a_0,b_0<1$, $0\le \alpha\le \pi/2$.}}
$$  %
Then %
$$ %
\varphi_r(z_1)=ae^{i\alpha},\quad \varphi_r(z_2)=be^{-i\alpha},
$$ %
where $a=a(r)$, $b=b(r)$ are functions of $r$.

Consider the function
$F(r)=\left(p_\mathbb{D}(ae^{i\alpha},be^{-i\alpha})\right)^2$.
After some
algebra, we find that %
\begin{equation}  \label{Equation 3.46}%
F(r)=\frac{a^2(r)-2a(r)b(r)\cos\alpha+b^2(r)}{(1-a(r)b(r))^2}.
\end{equation} %
Since for fixed $s$ and $\theta$, $r_0\le s<1$, $\theta\in
\mathbb{R}$, the hyperbolic distance between the points
$re^{i\theta}$ and $\varphi_r(se^{i\theta})$ is the same for all
$r$ in the interval $0\le r\le r_0$, it follows that %
\begin{equation} \label{Equation 3.47} %
\frac{da}{dr}=\frac{1-a^2}{1-r^2} \quad {\mbox{and}} \quad
\frac{db}{dr}=\frac{1-b^2}{1-r^2}.
\end{equation}  %
Differentiating (\ref{Equation 3.46}) and using formulas
(\ref{Equation 3.47}), we find:
\begin{equation}  \label{Equation 3.48}%
\frac{dF}{dr}=\frac{2}{(1-r^2)(1-ab)^3}\,\Phi(a,b,\alpha),
\end{equation}  %
where $\Phi=\Phi(a,b,\alpha)$ is the following function:%
$$ %
\Phi=(a(1-a^2)-(a(1-b^2)+b(1-a^2))\cos\alpha+b(1-b^2))(1-ab)
$$
$$ %
\hspace{-1.5cm}+(a^2-2ab\cos\alpha+b^2)(a(1-b^2)+b(1-a^2)). \ \ \
\
$$ %

It is clear that $\Phi(a,b,\alpha)$ is an increasing function of
$\alpha$ on the interval $0\le \alpha\le \pi/2$. After simple
calculation, we find that $\Phi(a,b,0)=0$ and therefore %
$$ %
\Phi(a,b,\alpha)\ge \Phi(a,b,0)=0.
$$ %
This, together with (\ref{Equation 3.48}), implies that
$\frac{dF(r)}{dr}\ge 0$ and therefore the pseudo-hyperbolic
distance in equation (\ref{Equation 3.45}) decreases, when $r$
decreases from $r_0$ to $0$. Therefore, our claim that $\varphi_r$
is a hyperbolic contraction is proved.  By
Proposition~\ref{Proposition 2.19}, the latter implies that ${\rm
cap}(E(r))$ decreases when $r$ decreases
from $r_0$ to $0$, which proves the theorem. %
\hfill  $\Box$

\section{Extremal properties of hedgehogs on 
evenly distributed radial intervals} 

In this section we consider problems with extremal configurations
lying on $n\ge 2$ radial intervals $I^*_k=\{z=te^{2\pi
i(k-1)/n}:\,0\le t\le 1\}$, $k=1,\ldots,n$. Since the intervals
$I_k$ are evenly distributed in $\mathbb{D}$ it is expected that
extremal configurations possess rotational symmetry by angle
$2\pi/n$.

First, we prove a theorem that generalizes the second inequality
of Theorem~\ref{Theorem 3.35} for sets lying on $2n\ge 4$
diameters.

\begin{theorem} \label{Theorem 4.1}%
Let $0<r_k<1$, $k=1,2$ and let $r$, $0<r<1$, be such that %
\begin{equation} \label{Equation 4.2}%
\uptau=\ell_{\mathbb{D}}([0,r])=\frac{1}{2}(\ell_{\mathbb{D}}([0,r_1])+\ell_{\mathbb{D}}([0,r_2]))%
 \end{equation} %
and, for $n\ge 2$, let %
$$ %
E_n(r_1,r_2)= (\cup_{k=0}^{n-1} e^{\pi
ik/n}[-r_1,r_1])\cup(\cup_{k=0}^{n-1} e^{\pi
i(2k+1)/2n}[-r_2,r_2]). %
$$ %
 Then %
\begin{equation} \label{Equation 4.3}%
{\rm
cap}(E_n(r_1,r_2))=8n\frac{\mathcal{K}(\kappa)}{\mathcal{K}'(\kappa)},
\quad {\mbox{where}} \quad
\kappa=\frac{r_1^{2n}+r_2^{2n}}{1+r_1^{2n}r_2^{2n}}.
\end{equation} %

 Furthermore, the following inequality holds:
\begin{equation} \label{Equation 4.4} %
{\rm cap}(E_n(r_1,r_2))\ge {\rm cap}(\cup_{k=0}^{2n-1} e^{\pi
ik/2n}[-r,r]). %
\end{equation}  %
Equality occurs in (\ref{Equation 4.4}) if and only if $r_1=r_2$.
\end{theorem}

\smallskip

\noindent %
\emph{Proof.} To establish (4.3), we use the function
    $\varphi=\varphi_2\circ \varphi_1$,  where $\varphi_1(z)=z^{2n}$,
    $\varphi_2(z)=(z+r_2^{2n})/(1+r_2^{2n}z)$, which  maps the sector $S=\{z\in \mathbb{D}:\,|\arg
    z|<\pi/2n\}$ conformally onto $\mathbb{D}$ slit along the interval
    $[-1,0]$. Using symmetry properties of $E_n(r_1,r_2)$, we find
    that %
   $$ 
{\rm cap}(E_n(r_1,r_2))=2n\,{\rm
cap}([0,(r_1^{2n}+r_2^{2n})/(1+r_1^{2n}r_2^{2n})]).
   $$ 
This together with (\ref{Equation 1.7}) gives (\ref{Equation
4.3}).

\smallskip

One way to prove the monotonicity property of ${\rm
cap}(E_n(r_1,r_2))$ is to differentiate the function in
(\ref{Equation 4.3}) and check if its derivative is positive.
Here, we demonstrate a different approach, which may be useful
when an explicit expression for the derivative is not known.

The proof presented below is similar to the proof of the
second inequality in (\ref{Equation 3.36}) in Theorem~\ref{Theorem 3.35}. We use the function %
$$ %
g(z)=\frac{z^{2n}+r_2^{2n}}{1+r_2^{2n}z^{2n}}
$$ %
to map the domain $\{z\in \mathbb{D}\setminus E_n(r_1,r_2):\,|\arg
z|<\pi/2n\}$ conformally onto the unit disk slit along the
interval
$(-1,F_n]$, where $F_n=F_n(r_1,r_2)$ is defined as %
$$ %
F_n=\frac{r_1^{2n}+r_2^{2n}}{1+r_1^{2n}r_2^{2n}}.
$$ %

Since $E_n(r_1,r_2)$ possesses $2n$-fold rotational symmetry about
$0$, it follows from the symmetry principle for the
module of family of curves, that %
$$ %
{\rm cap}(E_n(r_1,r_2))=2n\,{\rm cap}([0,F_n(r_1,r_2]).
$$ %
Therefore, to minimize ${\rm cap}(E_n(r_1,r_2))$ under the
constraint (\ref{Equation 4.2}), we can minimize $F_n$ under the
same constraint.

Using the variables $u=(1+r_1)/(1-r_1)$, $v=(1+r_2)/(1-r_2)$,
constrained by the condition $uv=C$, we express $F_n$ as follows: %
$$ %
F_n=\frac{(C-1+(u-v))^{2n}+(C-1+(v-u))^{2n}}{(C+1+(u+v))^{2n}+(C+1-(u+v))^{2n}}.
$$ %

To minimize $F_n(u,v)$ under the constraint $uv=C$, we introduce
Lagrange's function %
$$  %
L=F_n(u,v)+\lambda uv,\quad \lambda\in \mathbb{R}.
$$  %
Differentiating this function, we find %
\begin{equation}  \label{Equation 4.5}%
\frac{\partial L}{\partial u}-\lambda v= \frac{\partial
L}{\partial v}-\lambda u=M(u,v,C),
\end{equation}  %
where $M=M(u,v,C)$ is defined as follows:%
$$  %
M=2n\,\left((C+1+(u+v))^{2n}+(C+1-(u+v))^{2n}\right)^{-2}\,(M_1N_1-M_2N_2),
$$  %
$$ %
{\hspace{-0.85cm}}M_1=(C+1+(u+v))^{2n}+(C+1-(u+v))^{2n},
$$ %
$$  %
{\hspace{-0.85cm}}M_2=(C-1+(u-v))^{2n}+(C-1+(v-u))^{2n},
$$ %
$$  %
N_1=(C-1+(u-v))^{2n-1}-(C-1+(v-u))^{2n-1},
$$ %
$$ %
N_2=(C+1+(u+v))^{2n-1}-(C+1-(u+v))^{2n-1}.
$$ %

It follows from equation (\ref{Equation 4.5}) that if $(u,v)$ is a
critical point of the minimization problem under consideration,
then
$u=v=\sqrt{C}$. In this case  %
$$ 
F_n(\sqrt{C},\sqrt{C})=\frac{2(C-1)^{2n}}{(\sqrt{C}+1)^{4n}+(\sqrt{C}-1)^{4n}}.
$$ 
Since $u\ge 1$, $v\ge 1$, $uv=C$ and there is only
one critical point $(u,v)=(\sqrt{C},\sqrt{C})$ of the minimization
problem under consideration and since $F(u,v)=F(v,u)$, it follows
that $F(u,v)$ achieves its minimal value either at the point
$(u,v)=(\sqrt{C},\sqrt{C})$ or at the point $(u,v)=(1,C)$. In the
latter case, we have%
$$ 
F_n(1,C)=\frac{(C-1)^{2n}}{(C+1)^{2n}}. %
$$ 
The inequality $F_n(\sqrt{C},\sqrt{C})<F_n(1,C)$ is equivalent to
the
following inequality: %
$$ %
(\sqrt{C}+1)^{4n}+(\sqrt{C}-1)^{4n}>2(C+1)^{2n}.
$$ %
Using binomial expansion, this inequality can be written as %
$$ %
\sum_{k=0}^{2n} \binom{4n}{2k}C^{2n-k}> \sum_{k=0}^{2n}
\binom{2n}{k}C^{2n-k}.
$$ %
Since $\binom{2n}{2k}>\binom{n}{k}$ for all $n\ge 2$ and $1\le
k\le n-1$, the latter inequality holds true.

Thus, $F_n(u,v)$ takes its minimal value when $u=v=\sqrt{C}$ and
therefore the inequality (\ref{Equation 4.4}) is proved.
 \hfill  $\Box$


\begin{rem} %
The binomial inequalities similar to the one used in the proof of
Theorem~\ref{Theorem 4.1} are known to the experts. To prove that
$\binom{2n}{2k}>\binom{n}{k}$, we can argue as follows.

Let $n\ge 2$ and $0\le k\le n-k$. If $k=0$, then
$\binom{2n}{2k}=\binom{n}{k}=1$. Suppose that $\binom{2n}{2k}\ge
\binom{n}{k}$ for $0\le k\le n-k-1$. Then %
$$ %
\binom{2n}{2(k+1)}=\binom{2n}{2k}\frac{(2n-2k)(2n-2k-1)}{(2k+1)(2k+2)}\ge
\binom{n}{k}\frac{2(n-k)(2n-2k-1)}{2(k+1)(2k+1)}  %
$$ %
$$ %
\hspace{-2.5cm}=\binom{n}{k+1}\frac{2n-2k-1}{2k-1}>\binom{n}{k+1}. %
$$ %
Now the required inequality follows by induction.
\end{rem}

\smallskip

It was shown in the proof of Theorem~\ref{Theorem 4.1} that there
is only one critical point in the minimization problem considered
in that theorem. Therefore, the following monotonicity result is
also proved.

 \begin{corollary} \label{Corollary 4.10}  %
Under the assumptions of Theorem~\ref{Theorem 4.1}, suppose that
$r_1=s$ and $r_2=r_2(s)$ is such that condition (\ref{Equation
4.2}) holds. Then ${\rm cap}(E_n(s,r_2(s)))$ strictly decreases
from $8n\frac{\mathcal{K}(\kappa_0)}{\mathcal{K}'(\kappa_0)}$ to
$16n\frac{\mathcal{K}(\kappa_1)}{\mathcal{K}'(\kappa_1)}$, when
$s$ runs from $0$ to
$r$, where %
$$  %
\kappa_0=\left(\frac{e^{2\uptau}-1}{e^{2\uptau}+1}\right)^{2n},
\quad \quad
\kappa_1=\left(\frac{e^{\uptau}-1}{e^{\uptau}+1}\right)^{4n}.
$$  %
 \end{corollary} %

\medskip

Since the conformal capacity is conformally invariant, the
conformal capacity of an interval $E$ on $(-1,1)$ remains constant
when $E$ moves along the diameter $(-1,1)$ so that its hyperbolic
length is fixed. For $n$ intervals situated on $n$ equally
distributed radial intervals the latter property is not true any
more but, as our next theorem shows, if all these intervals have
equal hyperbolic lengths and move synchronically, the conformal
capacity of their union changes monotonically. Actually, the
non-strict monotonicity is already established in
Theorem~\ref{Theorem 3.44}. Thus, our intention here is to prove
the strict monotonicity result and relate it with certain
properties of relevant transcendental functions.

\begin{theorem} \label{Theorem 4.8}%
        Let $E_\uptau(r)$ be a closed subinterval of $[0,1)$ with the initial point at $0<r<1$ and hyperbolic length
        $\uptau>0$. For $n\in \mathbb{N}$, let $E_\uptau^n(r)= \cup_{k=1}^n
e^{2\pi i(k-1)/n}E_\uptau(r)$.  Then the conformal capacity ${\rm cap}(E_\uptau^n(r))$ is given by%
\begin{equation} \label{Equation 4.9}%
{\rm cap}(E_\uptau^n(r))=4n
\frac{\mathcal{K}(\kappa)}{\mathcal{K}'(\kappa)},
\end{equation}  %
where %
\begin{equation} \label{Equation 4.10}
\kappa=\frac{\rho^n-r^n}{1-r^n\rho^n},\quad \quad
\rho=\frac{e^{\uptau}(1+r)-(1-r)}{e^{\uptau}(1+r)+(1-r)}.
\end{equation} 

Furthermore, ${\rm cap}(E_\uptau^n(r))$ strictly increases from
$4n\frac{\mathcal{K}(\kappa_0)}{\mathcal{K}'(\kappa_0)}$ with
$\kappa_0=(e^{\uptau}-1)^n/(e^{\uptau}+1)^n$  to
$4n\frac{\mathcal{K}(\kappa_1)}{\mathcal{K}'(\kappa_1)}$ with
$\kappa_1=(e^{\uptau}-1)/(e^{\uptau}+1)$, when $r$ varies from $0$
to $1$.
  \end{theorem}


    \noindent %
    \emph{Proof.} As in the proof of Theorem~\ref{Theorem 4.1}, we use the function
    $\varphi=\varphi_2\circ \varphi_1$,  where $\varphi_1(z)=z^n$,
    $\varphi_2(z)=(z-r^n)/(1-r^nz)$.  Then $\varphi$ maps the sector $S=\{z\in \mathbb{D}:\,|\arg
    z|<\pi/n\}$ conformally onto $\mathbb{D}$ slit along the interval
    $[-1,0]$. Furthermore, $\varphi$ maps $E_\uptau(r)$ onto the interval $[0,\rho]$ with $\rho$ defined as in
    (\ref{Equation 4.10}). Using symmetry properties of $E_\uptau^n(r)$, we find
    that %
 $$ 
{\rm cap}([r^n,\rho^n])={\rm cap}([0,(\rho^n-r^n)/(1-r^n\rho^n)]).
$$ 
This equation together with (\ref{Equation 1.7}) gives
(\ref{Equation 4.9}).

\smallskip

As we have mentioned in the proof of Theorem~\ref{Theorem 4.1}, we
know two approaches to prove the monotonicity property of the
conformal capacity in that theorem. The same approaches can be
used to prove the monotonicity statement of Theorem~\ref{Theorem
4.8}. Here, we demonstrate one more approach, which also may be
useful when an explicit expression for the derivative is not
known. We first note that ${\rm cap}(E_\uptau^n(r))$ is an
analytic function of $r$. This follows from equation
(\ref{Equation 4.9}). Since ${\rm cap}(E_\uptau^n(r))$ is not
constant and analytic, it follows that ${\rm cap}(E_\uptau^n(r))$
is not constant on any subinterval of $[0,1)$. Furthermore, it
follows from Theorem~\ref{Theorem 3.44} that ${\rm
cap}(E_\uptau^n(r))$ is a non-decreasing function. Since it is
non-decreasing and not
constant on any interval, it is strictly increasing.  %
\hfill $\Box$

\smallskip

Our next theorem can be considered as a counterpart of the
subadditivity property of the conformal capacity discussed in
Proposition~\ref{Proposition 2.3.1}.


\begin{theorem} \label{Theorem 4.11}%
        Let $E_k$, $k=1,\ldots,n$, be compact sets on the interval $I=[0,1)$ having positive logarithmic capacities
        and such that every point of $E_k$ is regular for the  Dirichlet problem
in $\mathbb D\setminus E_k$, $k=1,2,\dots,n$. Then
        \begin{equation}\label{Equation 4.12}
       \frac{1}{n}\sum_{k=1}^n{\rm cap}(\cup_{j=1}^n e^{2\pi i(j-1)/n}E_k)\le
       {\rm cap}(\cup_{k=1}^n e^{2\pi i(k-1)/n}E_k)<\sum_{k=1}^n{\rm cap}
        (E_k).
        \end{equation} %

        Equality occurs in the first inequality if and only if for each $k$ and $j$, $E_k$ coincides
        with $E_j$ up to a set of zero logarithmic capacity.
    \end{theorem}

    \smallskip

    \noindent %
    \emph{Proof.} The second inequality is just  the
    subadditivity property of the conformal capacity stated in Proposition~\ref{Proposition 2.3.1}.

    To prove the
    first inequality, we use the method of separation of
    components of a condenser in the style of Dubinin's paper
    \cite{D1}.
    Let $u$ denote the potential function of the condenser $(\mathbb{D},\bigcup_{k=1}^n e^{2\pi
    i(k-1)/n}E_k)$. Since every point of $E_k$ is regular for the  Dirichlet problem, it follows that
    $u$ is continuous on $\overline{\mathbb{D}}$. Let $u'_k$ and $u''_k$ denote the functions
    obtained from $u$, first by restricting $u$, respectively, onto the sector
    $S_1=\{z\in \mathbb{D}:\,\pi (2k-3)/n\le \arg z\le 2\pi(k-1)/n\}$
    or onto the sector $S_2=\{z\in \mathbb{D}:\,2\pi (k-1)/n\le \arg z\le \pi(2k-1)/n\}$
    and then extending this restriction by symmetry on the whole
    unit disk. Then each of the functions $u'_k$ and $u''_k$ is admissible for the
    condenser $(\mathbb{D},\cup_{j=1}^n e^{2\pi i(j-1)/n}E_k)$.
    Furthermore, each of these functions possesses $n$-fold
    rotational symmetry about $0$ and is symmetric with respect to
    the real axis. Therefore, the following inequality holds:
    \begin{equation} \label{Equation 4.13} %
\frac{1}{n}{\rm cap}(\cup_{j=1}^n e^{2\pi i(j-1)/n}E_k)\le
\int_{S_1\cup S_2} |\nabla u_k|^2\,dm, \quad k=1,\ldots,n.
    \end{equation}  %
Summing up all the inequalities in (\ref{Equation 4.13}), we
obtain the first inequality in (\ref{Equation 4.12}).

    Furthermore, since every point of the sets $E_j$, $j=1,\ldots,n$, is
regular for the Dirichlet problem, it follows that $u'_k$ or
$u''_k$ defined above in the proof is a
    potential function of $(\mathbb{D},\cup_{j=1}^n e^{2\pi
    i(j-1)/n}E_k)$ if and only of $E_j=E_k$ for all
    $j=1,\ldots,n$. Therefore, if $E_j\not = E_k$ for some
    $j\not=k$, then we have the strict inequality in (\ref{Equation 4.13}) and in
    the first inequality in (\ref{Equation 4.12}) as well.
\hfill $\Box$

\medskip

Above we discussed results on the conformal capacity of compact
sets lying on a fixed number of radial intervals. In our next
theorem, we work with compact sets on $m\ge n\ge 2$ radial
intervals that are ``densely spread'' over $\mathbb{D}$ in the
sense that the angle between any two neighboring intervals is $\le
2\pi/n$.


\begin{theorem} \label{Theorem 4.14}%
Let $E_0=\overline{\mathbb{D}}_r$, $0<r<1$, and let $E_1\subset
[r,1)$ be a compact set of hyperbolic length $\uptau>0$. Let
$0=\alpha_0=\alpha_1<\alpha_2<\ldots<\alpha_m<\alpha_{m+1}=2\pi$
be such that $\alpha_{k+1}-\alpha_k\le 2\pi/n$, $k=1,\ldots,m$,
$m\ge n\ge 2$. Consider compact sets $E=\cup_{k=0}^m
e^{i\alpha_k}E_1$, $E^*=\cup_{k=0}^n e^{i 2\pi(k-1)/n}E_1$, and
$E_n(r,\uptau)=E_0\cup(\cup_{k=1}^n e^{i 2\pi(k-1)/n}[r,r_1])$,
where $r_1$ is such that $r<r_1<1$ and
$\ell_{\mathbb{D}}([r,r_1])=\uptau$.
Then %
\begin{equation} \label{Equation 4.15}%
{\rm cap}(E)\ge {\rm cap}(E^*)\ge {\rm cap}(E_n(r,\uptau)).
\end{equation}  %

Equality in the first inequality occurs if and only if $m=n$.
    \end{theorem}

    \smallskip

    \noindent %
    \emph{Proof.}  The first inequality in (\ref{Equation 4.15}), together with the statement on the equality cases,
     follows from Theorem~5 in \cite{BS}. Then the second inequality
  follows from the contraction principle stated in Theorem~\ref{Theorem 3.44}.   \hfill  $\Box$

\medskip



    Polarization and the contraction principle are customarily
    applied when lower bounds for the conformal capacity are
    needed. In our next theorem, we present a result with an
    upper bound for this capacity.

\begin{theorem} \label{Theorem 4.19}%
For $0\le r\le a<b<1$, let $E_0=\overline{\mathbb{D}}_r$ and let
$E_1\subset [a,b]$ be a compact set of hyperbolic length
$\uptau>0$, $0<\uptau<\ell_{\mathbb{D}}([a,b])$. Let
$0=\alpha_1<\alpha_2<\ldots<\alpha_n<\alpha_{n+1}=2\pi$. Let %
$$  %
E=E_0\cup (\cup_{k=1}^n e^{i\alpha_k}E_1),\quad E^*=E_0\cup
(\cup_{k=1}^n e^{i2\pi(k-1)/n}E_1), \quad E^{a,b}=E_0\cup
(\cup_{k=1}^n e^{i2\pi(k-1)/n}[a,b]).
$$  %
Then %
\begin{equation} \label{Equation 4.17}%
{\rm cap}(E)\le {\rm cap}(E^*)< {\rm cap}(E^{a,b}).
\end{equation}  %

Equality in the first inequality in (\ref{Equation 4.17}) occurs
if and only if $\alpha_k=2\pi(k-1)/n$, $k=1,\ldots,n$.
    \end{theorem} %

    \smallskip

    \noindent %
    \emph{Proof.}  The first inequality in (\ref{Equation 4.17}) together with the statement on the cases of equality
    follows from Dubinin's dissymmetrization  results; see Theorem~4.14 in \cite{Du}.
    Then, since the conformal capacity is an increasing function of a set,  the second inequality
    follows.%
  \hfill  $\Box$


\medskip

\begin{rem} %
If $E_1$ in Theorem~\ref{Theorem 4.19} is an interval
$[c,d]\subset [a,b]$, then the upper bound in (\ref{Equation
4.17}) can be replaced with ${\rm cap}(E_0\cup (\cup_{k=1}^n
e^{i2\pi(k-1)/n}[s,b]))$ with $s$ in $(a,b)$ such that
$\ell_{\mathbb{D}}([s,b])=\uptau$. This follows, for instance, from the monotonicity property stated in Theorem~\ref{Theorem 3.44}. %
\end{rem}

  \section{Symmetrization transformations in hyperbolic metric}

In this section, we discuss possible counterparts of classical
symmetrization-type transformations applied with respect to the
hyperbolic metric. We note that Steiner, Schwarz and circular
symmetrizations destroy hedgehog structures, in general,   and
therefore their applications to problems studied in previous
sections of this paper are limited. We define symmetrizations
using the following notations, which are convenient for our
purposes:
$L_\alpha=\{z:\,\Im(e^{-i\alpha}z)=0\}$, $L^a=\{z:\,\Im z=a\}$,
$C_r=\{z:\,|z|=r\}$, $R_\alpha=\{z=te^{i\alpha}:\,t\ge0\}$. First,
we mention results obtained with  Steiner symmetrization.
\begin{definition} \label{Definition 5.1}%
        Let $E\subset \mathbb{C}$ be a compact set. The Steiner
        symmetrization of $E$ with respect to the imaginary axis
        is defined to be the compact set
        \begin{equation}\nonumber
        E^*=\{z=x+iy:\,E\cap
        L^y\not=\emptyset, |x|\le
        (1/2)\ell(E\cap L^y)\},
        \end{equation}
         where $\ell(\cdot)$ stands for the one-dimensional Lebesgue measure. %

        Furthermore, the Steiner symmetrization of $E$ with respect to the
        line $L_\alpha$ is defined to be the compact set
        $E_\alpha^*=e^{i(\alpha-\pi/2)}\left(e^{i(\pi/2-\alpha)}E\right)^*$. %
        \end{definition}

For the properties and results obtained with the Steiner
symmetrization, the interested reader may consult \cite{PS},
\cite{Du}, \cite{B}.

We note that if $E$ is a compact set in $\mathbb{D}$, then
$E_\alpha^*\subset \mathbb{D}$ for all $\alpha\in \mathbb{R}$ and
therefore ${\rm cap}(E_\alpha^*)$ is well defined. As is well
known, Steiner symmetrization does not increase the capacity of a
condenser and therefore it does not increase the conformal
capacity. Also, Steiner symmetrization preserves Euclidean area
but, in general,  it strictly decreases the hyperbolic area
$A_\mathbb{D}(E)$ of $E$ that is defined by (\ref{Hyperbolic
area}). Therefore, it is not an equimeasurable rearrangement with
respect to the hyperbolic metric. Thus, to study problems on the
hyperbolic plane, a version of Steiner symmetrization, which
preserves the hyperbolic area and does not increase the conformal
capacity, is needed.
To define this symmetrization, we will use the function %
$$ 
 w=\varphi_0(z)=\log \frac{1+z}{1-z}, \quad \varphi_0(0)=0, %
 $$ 
 which maps $\mathbb{D}$ conformally onto the horizontal strip $\Pi=\{w:\,|\Im w|<\pi/2\}$.
 We note that $\varphi$ maps the hyperbolic geodesic $(-1,1)$
 onto the real axis and it maps the curves equidistant from $(-1,1)$,
 (which are circular arcs in $\mathbb{D}$ joining the points $1$ and $-1$),
 onto the horizontal lines $\{w:\,\Im w=c\}$, $0<|c|<\pi/2$, in $\Pi$.
\begin{definition} \label{Definition 5.3} %
        Let $E$ be a compact set in $\mathbb{D}$. The hyperbolic Steiner
        symmetrization of $E$ with respect to the hyperbolic geodesic
        $(-i,i)$, centered at $z=0$, is defined as
      $$ 
        E_h^*=\varphi_0^{-1}\left(\left(\varphi_0(E)\right)^*\right),%
        $$ 
        where $\left(\varphi_0(E)\right)^*$ stands for the Steiner
        symmetrization as in Definition~\ref{Definition 5.1}.

        Furthermore, the hyperbolic Steiner symmetrization of $E$ with
        respect to a hyperbolic geodesic $\gamma$, centered at the point
        $a\in \gamma$, is defined as
        $$ %
        E_h(\gamma,a)=\psi^{-1}\left(\left(\psi(E)\right)_h^*\right),
        $$
        where $\psi(z)$ denotes the M\"{o}bius automorphism of
        $\mathbb{D}$, which maps $\gamma$ onto the hyperbolic geodesic
        $(-i,i)$ such that $\psi(a)=0$. %
        \end{definition}

The hyperbolic Steiner symmetrization was introduced by A.~Dinghas
\cite{Dinghas}. We note that this transformation symmetrizes sets
along the hyperbolic equidistant lines and not along the
hyperbolic geodesics. Previously it was used in several research
papers, see, for instance, \cite{KarpPeyerimhoff2002},
\cite{Gueritaud}. In our next theorem, we collect properties of
the hyperbolic Steiner symmetrization that are relevant to our
study.

\begin{theorem}[{see \cite{KarpPeyerimhoff2002}}] \label{Theorem 5.5}
        Let $E_h(\gamma,a)$ be the image of a compact set $E\subset
        \mathbb{D}$ under the hyperbolic Steiner symmetrization with
        respect to a hyperbolic geodesic $\gamma$, centered at $a\in
        \gamma$. Let $\gamma_\perp\ni a$ be the hyperbolic geodesic
        orthogonal to $\gamma$. Then the following hold true:
        \begin{enumerate} %
                \item[{\rm (1)}] If $E$ is a hyperbolic disk, then $E_h(\gamma,a)$ is
                a hyperbolic disk of the same hyperbolic area as $E$ having its center on $\gamma$. %
                \item[{\rm (2)}] $E_h(\gamma,a)$ is a compact set that is symmetric with respect to $\gamma$. %
                \item[{\rm (3)}] $\ell_\mathbb{D}(E_h(\gamma,a)\cap
                \gamma_\perp)=\ell_\mathbb{D}(E\cap
                \gamma_\perp)$. %
                \item[{\rm (4)}] $A_\mathbb{D}(E_h(\gamma,a))=A_\mathbb{D}(E)$. %
                \item[{\rm (5)}] ${\rm cap}(E_h(\gamma,a))\le {\rm cap}(E)$ with the
                sign of equality if and only if $E_h(\gamma,a)$ coincides with $E$
                up to a set of zero logarithmic capacity and up to a M\"obius automorphism of $\mathbb D$ that preserves $\gamma_\perp$. %
        \end{enumerate} %
\end{theorem}

\smallskip

\noindent %
\emph{Proof.} We use the notation that we set in the definitions
of Steiner and hyperbolic Steiner symmetrizations. We equip the
strip $\Pi$ with the hyperbolic metric $\lambda_\Pi(w)|dw|$
induced by the hyperbolic metric in $\mathbb D$ via the conformal
mapping $\varphi_0:\mathbb D\to\Pi$. That is, we have %
$$ %
\lambda_\Pi(w)|dw|=\lambda_{\mathbb
        D}(z)|dz|,\;\;\;w=\varphi_0(z),\;z\in\mathbb D,\;w\in\Pi.
$$ %
So, trivially, $\varphi_0$ is a hyperbolic isometry from $\mathbb
D$ to $\Pi$.

\smallskip

(1) Let $E$ be a hyperbolic disk in $\mathbb D$. Then $\psi(E)$ is
a hyperbolic disk in $\mathbb D$ and $\varphi_0\circ\psi(E)$ is a
hyperbolic disk in $\Pi$. It is easy to observe that hyperbolic
disks in $\Pi$ are horizontally convex (namely, their intersection
with any horizontal line is either empty or a single horizontal
rectilinear interval). It follows from the definition of Steiner
symmetrization that $(\varphi_0\circ\psi(E))^*$ is obtained by a
horizontal rigid motion of $\varphi_0\circ\psi(E)$ and it is a
hyperbolic disk in $\Pi$, symmetric with respect to the imaginary
axis. Since both $\psi^{-1}$ and $\varphi_0^{-1}$ are hyperbolic
isometries and preserve symmetries, $E_h(\gamma,a)$ is a
hyperbolic disk in $\mathbb D$, symmetric with respect to
$\gamma$. Moreover, a horizontal motion in $\Pi$ preserves the
hyperbolic area. Therefore, $A_{\mathbb
        D}(E_h(\gamma,a))=A_{\mathbb D}(E)$.

\smallskip

(2) It is well known that Steiner symmetrization transforms
compact sets to compact sets. So, if $E\subset \mathbb D$ is
compact, $E_h(\gamma,a)$ is compact, too. Moreover,
$(\varphi_0(E))^*$ is a set in $\Pi$, symmetric with respect to
the imaginary axis. Hence, $E_h^*$ is a compact set in $\mathbb
D$, symmetric with respect to the geodesic $(-i,i)$. It follows
that $E_h(\gamma,a)$ is symmetric with respect to $\gamma$.

\smallskip

(3) Since $\psi$ is a hyperbolic isometry on $\mathbb D$, the set
$\psi(E\cap \gamma_\perp)$ is a compact subset of the diameter
$(-1,1)$ having the same hyperbolic length as $E\cap\gamma_\perp$.
Therefore, $\varphi_0\circ\psi(E\cap\gamma_\perp)$ lies on the
real axis and
$\ell_\Pi(\varphi_0\circ\psi(E\cap\gamma_\perp))=\ell_{\mathbb
        D}(E\cap\gamma_\perp)$. The hyperbolic length (in $\Pi$) on the real
        axis is proportional to the Euclidean length $\ell$.
        Hence $\ell_\Pi(\varphi_0\circ\psi(E\cap\gamma_\perp))=\ell_\Pi((\varphi_0\circ\psi(E\cap\gamma_\perp))^*)$.
Thus the equality in (3) follows at once.

\smallskip

(4) This is true because $\varphi_0$, $\psi$, and Steiner
symmetrization (with respect to the imaginary axis, in $\Pi$)
preserve hyperbolic areas.

\smallskip

(5) The conformal maps $\varphi_0$ and $\psi$ preserve the
capacity of condensers. So, to prove the inequality in (5), it
suffices to show that for every compact subset $F$ of $\Pi$, we
have ${\rm cap}(\Pi,F)\geq {\rm cap}(\Pi, F^*)$. This is a
well-known symmetrization theorem \cite[Theorem 4.1]{Du}. The
equality statement follows from \cite[Theorem 1]{BP2012}. \hfill
$\Box$

\smallskip

\begin{rem} \label{Remark 5.6}%
        Of course, any transformation bearing the name ``symmetrization''
        must possess property (1) of Theorem~\ref{Theorem 5.5}; otherwise it is not a
        symmetrization. On the other side, it was shown by L.~Karp
        and N.~Peyerimhoff in \cite{KarpPeyerimhoff2002} and by
        F.~Gu\'{e}ritaud in \cite{Gueritaud} that even in the best case
        scenario the hyperbolic Steiner symmetrization changes hyperbolic triangles into
        sets that are not convex with respect to hyperbolic metric and
        therefore these sets are not hyperbolic triangles. In fact, the
        image $T_h^*$ of the hyperbolic triangle $T$ with vertices
        $z_1=-r$, $z_2=r$, $0<r<1$, and $z_3=-is+\sqrt{1+s^2}e^{i\beta}$,
        $s>0$, $\arctan s<\beta<\pi/2$, is a proper subset of the
        hyperbolic isosceles triangle $T_0$  with vertices $z_1=-r$,
        $z_2=r$, and $z_3=(\sqrt{1+s^2}-s)i$.
\end{rem} %

\smallskip

In relation to Theorem~\ref{Theorem 3.42} we suggest the following
problem, where two hyperbolic Steiner symmetrizations provide some
qualitative information about extremal configuration but these are
not enough to give a complete solution of the problem. This
problem is a counterpart of the problem for the logarithmic
capacity in the Euclidean plane, which was solved
in~\cite{BarnardPearceSolynin2002}.

\begin{problem} \label{Problem 5.7}%
        Find the minimal conformal capacity among all compact sets
        $E\subset \mathbb{D}$ with prescribed hyperbolic diameter $d>0$
        and prescribed hyperbolic area $0<A<4\pi \sinh^2(d/4)$.
        Describe possible extremal configurations.
\end{problem}


\medskip

The Schwarz symmetrization of $E$ in the plane with respect to a
point $a$ replaces $E$ by the disk centered at $a$ of the same
area. The hyperbolic analog of this symmetrization is the
following.

\begin{definition} \label{Definition 5.8}%
        Let $E$ be a compact set in $\mathbb{D}$. Then its hyperbolic
        Schwarz symmetrization with respect to $a\in \mathbb{D}$ is the
        hyperbolic disk, we call it $E_a^\#$, centered at  $a$ and such
        that $A_\mathbb{D}(E_a^\#)=A_\mathbb{D}(E)$. %
        \end{definition}

The hyperbolic Schwarz symmetrization was first suggested by
F.~Gehring \cite{G} and later used in \cite{FryntovRossi2001} and
\cite{Betsakos2013}. In particular, the following result was
proved.

\begin{theorem}[see \cite{G},\cite{FryntovRossi2001},\cite{Betsakos2013}] %
        If $E_a^\#$ is the hyperbolic Schwarz symmetrization of a compact set $E\subset \mathbb{D}$, then %
        $$ 
        {\rm cap}(E_a^\#)\le {\rm cap}(E)
        $$ 
        with the sign of equality if and only if  $E_a^\#$ coincides with
        $E$
        up to a set of zero logarithmic capacity and up to a M\"{o}bius automorphism of $\mathbb{D}$. %
        \end{theorem} %

\medskip

Next, we will discuss the hyperbolic circular  symmetrization of
$E$ with respect to the hyperbolic ray $[a,e^{i\alpha})_h$, $a\in
\mathbb{D}$, $\alpha\in \mathbb{R}$. This is defined as the image
of the interval $[0,1)$ under a M\"{o}bius automorphism $\varphi$
of $\mathbb{D}$ such that $\varphi(0)=a$,
$\varphi(1)=e^{i\alpha}$.

\begin{definition} \label{Definition 5.11} %
Let $E$ be a compact set in $\mathbb{D}$. Then its hyperbolic
        circular symmetrization with respect to the hyperbolic ray
        $[0,1)_h$ is a compact set $E_h^\circ \subset \mathbb{D}$  such
        that: (a) $0\in E_h^\circ$ if and only if $0\in E$, (b) for
        $0<r<1$, $E_h^\circ\cap C_r=\emptyset$ if and only if $E\cap
        C_r=\emptyset$, (c) if, for $0<r<1$, $E\cap C_r\not=\emptyset$,
        then $E_h^\circ\cap C_r$ is a closed circular arc on $C_r$ (which
        may degenerate to a point or may be the whole circle $C_r$)
        centered at $z=r$ such that
        $\ell_\mathbb{D}(E_h^\circ \cap C_r)=\ell_\mathbb{D}(E\cap C_r)$.%

        Furthermore, the hyperbolic circular symmetrization of $E$ with
        respect to a hyperbolic geodesic ray $[a,e^{i\alpha})_h$ is
        defined as        %

        \begin{equation} \label{Equation 5.9}%
        E_h^\circ(a,\alpha)=\varphi^{-1}\left(\left(\varphi(E)\right)_h^\circ\right),
        \end{equation} %
        where $\varphi(z)$ denotes the M\"{o}bius automorphism of
        $\mathbb{D}$, which maps $[a,e^{i\alpha})_h$ onto the interval
        $[0,1)$. %
\end{definition}

\begin{rem} \label{Remark 5.13} %
        We immediately admit here that the hyperbolic circular
        symmetrization with respect to $[0,1)_h$ is exactly the classical
        circular symmetrization with respect to the positive real axis.
        The reason for this is the fact that the hyperbolic density
        $\lambda_{\mathbb D}(z)$ is constant on circles centered at $0$. Thus, this
        transformation will not provide any new information that is not
        available via classical symmetrization methods. However, the
        variant in formula (\ref{Equation 5.9}) can be used to obtain additional
        information that may be
        rather interesting. %
\end{rem} %


It follows from Definition~\ref{Definition 5.11} that the
hyperbolic area and conformal capacity ${\rm
cap}(E_h^\circ(a,\alpha))$ do not depend on $\alpha$. For the
Euclidean area of $E_h^\circ(a,\alpha)$ we have the following
result.

\begin{lemma} \label{Lemma 5.14}%
        Let $E\subset \mathbb{D}$ be a compact set of positive area, and
        let $A(\alpha)$
        denote the Euclidean area
        of $E_h^\circ(r,\alpha)$, $0<r<1$, $0\le \alpha\le \pi$,
        considered as a function of $\alpha$. If $E_h^o(r,\alpha)$ does
        not coincide with a hyperbolic disk centered at $r$ up to measure
        $0$, then $A(\alpha)$ is strictly increasing in $\alpha$ on the
        interval $0\le \alpha\le \pi$; otherwise $A(\alpha)$ is constant
        for $0\le \alpha \le \pi$.

\end{lemma} %

\noindent %
\emph{Proof.} 
Let $E_0$ denote the circular symmetrization of $\varphi(E)$ with
respect to $[0,1)_h$, with $\varphi$ defined as in
Definition~\ref{Definition 5.11}.
Then $E_h^o(r,\alpha)=\psi(E_0)$, where  %
$$ %
\psi(z)=(e^{i\beta}z+r)/(1+e^{i\beta}rz), %
$$ %
with $\beta=\beta(\alpha)$ chosen such that %
$$ %
e^{i\alpha}=\psi(1)=(e^{i\beta}+r)/(1+e^{i\beta}r). %
$$ %
The function $\beta(\alpha)$ strictly increases from $0$ to $\pi$, when $\alpha$ increases from $0$ to $\pi$.
Therefore, to prove monotonicity of a certain characteristic $F$ of $E_h^o(r,\alpha)$ as a function of $\alpha$,
we can consider $\alpha=\alpha(\beta)$ as a function of $\beta$ and treat $F$ as a function of $\beta$.
Thus, we will work with the function
$A_1(\beta)=A(\alpha(\beta))$.

The Euclidean area $A_1(\beta)$ can be found as follows:
\begin{equation} \label{Equation 5.12} 
A_1(\beta)=\int_{E_0}|\psi'(z)|^2\,dm=(1-r^2)^2\int_0^1\left(\int_{-\theta(r)}^{\theta(r)}\frac{d\theta}{|1+r\rho
        e^{i(\theta+\beta)}|^4}\right)\,\rho\,d\rho, %
        \end{equation}  
        where $0\le \theta(r)\le \pi$ is defined by the condition $\gamma(r)=E_0\cap C_r=\{re^{i\theta}:\,|\theta|\le \theta(r)\}$.

If $0<\theta(r)<\pi$, then $\gamma(r)$ is a non-degenerate proper arc on
$C_r$ that is centered at $z=r$. Using this symmetry and the monotonicity property of the function %
$$ %
g(t)=1+2r\rho \cos t+r^2\rho^2, \quad 0\le t \le \pi,
$$ %
we conclude that the inner integral in (\ref{Equation 5.12}), that is the integral %
$$  %
I(\beta)=\int_{-\theta(r)}^{\theta(r)}\frac{d\theta}{|1+r\rho
        e^{i(\theta+\beta)}|^4}=\int_{-\theta(r)+\beta}^{\theta(r)+\beta}\frac{dt}{(1+2r\rho
        \cos t+r^2\rho^2)^4},
$$  %
is a strictly increasing function of $\beta$, $0\le \beta\le \pi$.
Therefore, if $E_0$ differs by positive area from the disk
centered at $0$, which has the same area as $E_0$, then the double
integral in the right-hand side of (\ref{Equation 5.12}) strictly
increases on the interval $0\le \beta\le \pi$.   This proves that
the area $A_1(\beta)=A(\alpha)$ strictly increases on the interval
$0\le \alpha\le \pi$.     Of course, if $E_h^o(r,\alpha)$
coincides with a hyperbolic disk centered at $r$ up to measure
$0$, then $A(\alpha)$ is constant for $0\le \alpha \le \pi$.\hfill
$\Box$

\medskip

Now we add two results in the spirit of our theorems on the
conformal capacity of hedgehogs proved in previous sections.
Figure~5 illustrates proofs of these results, which are given in
Theorems~\ref{Theorem 5.16} and \ref{Theorem 5.20} below. Examples
of an admissible set and an extremal set of Theorem~\ref{Theorem
5.16} are shown in parts (a) and (b) of Figure~5 and example of an
admissible and extremal sets of Theorem~\ref{Theorem 5.20} are
shown in parts (c) and (d) of this figure.

\begin{theorem} \label{Theorem 5.16}%
        For $0<r<1$ and $0\le \alpha\le \pi$, let
        $C_r(\alpha)=\{re^{i\theta}:\,|\theta|\le \alpha\}$ and let
        $E\subset \{z:\,r\le |z|<1\}$ be a compact set such that the
        length of its radial projection $E_{pr}(r)$ on the circle $C_r$ is
        $\ge
        2\alpha r$. Then %
        \begin{equation}  \label{Equation 5.14}%
        {\rm cap}(E)\ge {\rm cap}(C_r(\alpha))
        \end{equation} %
        with the sign of equality if and only if $E$ coincides with
        $C_r(\alpha)$ up to a set of zero logarithmic capacity and up to
        a rotation about $0$.
\end{theorem} %

\noindent %
\emph{Proof.}  Consider a function $\varphi(z)=rz/|z|$.
It can be easily shown that $\varphi$ is a hyperbolic contraction from $E$ onto $E_{pr}(r)$.
Hence, by Proposition~\ref{Proposition 2.19}, %
\begin{equation}  \label{Equation 5.15}%
{\rm cap}(E)\ge {\rm cap}(E_{pr}(r)).
\end{equation}  %
As well known, the circular symmetrization decreases the conformal capacity,
precisely, the following holds: %
\begin{equation}  \label{Equation 5.16}%
{\rm cap}(E_{pr}(r))\ge {\rm cap}(C_r(\alpha)) %
\end{equation}  %
with the sign of equality if and only if $E_{pr}(r)$ coincides
with $C_r(\alpha)$ up to a set of zero logarithmic capacity and a
rotation about $0$. Combining (\ref{Equation 5.15}) and
(\ref{Equation 5.16}), we obtain (\ref{Equation 5.14})
 together with the statement about cases of equality in it. \hfill     $\Box$ %

\smallskip

\begin{figure} 
\begin{minipage}{1.0\textwidth}
$$\includegraphics[scale=0.6,angle=0]{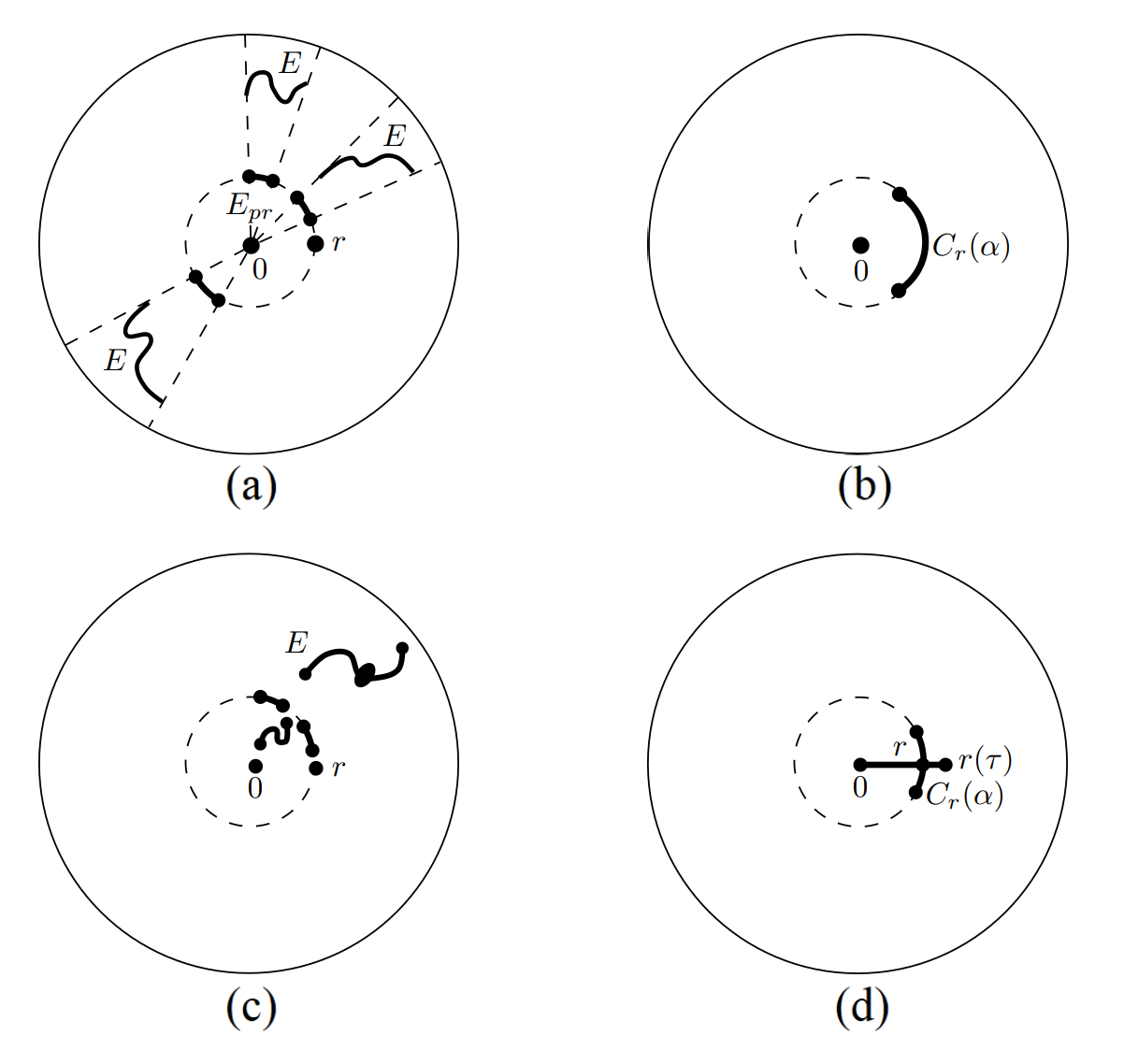}
$$
\end{minipage}
\caption{Admissible and extremal sets of Theorems~\ref{Theorem 5.16} and \ref{Theorem 5.20}.}%
\end{figure}

\begin{theorem} \label{Theorem 5.20}%
        For $0<r<1$, $0\le \alpha\le \pi$, and $\uptau\ge \uptau(r)=
        \log((1+r)/(1-r))$, let $E\subset \mathbb{D}$ be a compact set
        such that $\ell_{\mathbb{D}}(E\cap C_r)=2 \alpha r/(1-r^2)$ and
        $\ell_\mathbb{D}(E_{pr}^o)=\uptau$, where $E_{pr}^o$
        stands for the circular projection of $E$ onto the radius $[0,1)$.  Then %
        \begin{equation}  \label{Equation 5.18}%
        {\rm cap}(E)\ge {\rm cap}(C_r(\alpha)\cup [0,r(\uptau)]), \quad
        \quad {\mbox{where $r(\uptau)=(e^{\uptau}-1)/(e^{\uptau}+1)$,}}
        \end{equation} %
        with the sign of equality if and only if $E$ coincides with
        $C_r(\alpha)\cup [0,r(\uptau)]$ up to a set of zero logarithmic
        capacity and up to a rotation about $0$.
\end{theorem} %

\noindent %
\emph{Proof.} Performing the circular symmetrization as in the
proof of Theorem~\ref{Theorem 5.16}, we obtain the following
inequality:
\begin{equation} \label{Equation 5.19}%
{\rm cap}(E)\ge {\rm cap}(C_r(\alpha)\cup E_{pr}^o) %
\end{equation}  %
with the sign of equality if and only if $E$ coincides with
$C_r(\alpha)\cup E_{pr}^o$ up to a set of zero logarithmic
capacity and up to a rotation about $0$.

If $E_{pr}^o$ does not coincide with the interval $[0,r(\uptau)]$
up to a set of zero logarithmic capacity, then we can perform
polarization transformations as in the proof of Lemma~\ref{Lemma
3.4} to transform $E_{pr}^o$ into the interval $[0,r(\uptau)]$.
Furthermore, our assumption that
$\ell_{\mathbb{D}}(E_{pr}^o)=\uptau\ge \log((1+r)/(1-r))$
guarantees that these polarizations do not change $C_\alpha(r)$.
Hence, we obtain the inequality %
\begin{equation} \label{Equation 5.20}%
{\rm cap}(C_r(\alpha)\cup E_{pr}^o)\ge {\rm cap}(C_r(\alpha)\cup
[0,r(\uptau)])
\end{equation}  %
with the sign of equality if and only if $E_{pr}^o$ coincides with
$[0,r(\uptau)]$ up to a set of zero logarithmic capacity.
Combining (\ref{Equation 5.19}) and (\ref{Equation 5.20}), we
obtain (\ref{Equation 5.18}) together with the statement about
cases of equality in it. \hfill $\Box$

\smallskip

The restriction $\uptau\ge \uptau(r)$ in the previous theorem is
due limitations of  methods used in this paper. In relation with
this theorem, we suggest two problems.

\begin{problem} \label{Problem 5.24}%
        Find $\min_{E} {\rm cap}(E)$ over all compact sets $E\subset
        \mathbb{D}$ satisfying the assumptions of Theorem~\ref{Theorem 5.20} with
        $0<\uptau<\uptau(r)$. %
\end{problem}  %

\begin{problem} \label{Problem 5.25}%
        Let $E\subset \{z:\,r_1\le |z|\le r_2\}$, $0<r_1<r_2<1$, be a
        compact set such that its radial projection on the unit circle
        $\mathbb{T}$ coincides with $\mathbb{T}$ and its circular
        projection onto the radius $[0,1)$ coincides with the interval
        $[r_1,r_2]$, as it is shown in Figure~6, and let $E^*=\overline{\mathbb{D}}_{r_1}\cup
        [r_1,r_2])$. Is it true that ${\rm cap}(E)\ge {\rm cap}(E^*)$? If this is not true then identify the shape of $E$
        minimizing the conformal capacity under the assumptions of this
        problem. %
\end{problem}  %

\begin{rem} %
Let $\gamma=\gamma(r_1,r_2)\ni r_2$ denote the hyperbolic geodesic
tangent to the circle $C_{r_1}$. We assume that $\gamma$ touches
$C_{r_1}$ at the point $z_1=r_1e^{i\theta_0}$, $0<\theta_0<\pi/2$.
Let $\alpha$ denote the arc of $\gamma$ with the endpoints $z_1$,
$r_2$ and let $\overline{\alpha}=\{z:\,\bar{z}\in \alpha\}$. It
can be shown with the help of polarization that if $E^*$ minimizes
the conformal capacity in Problem~\ref{Problem 5.25}, then $E^*$
belongs to the compact set bounded by the arcs $\alpha$,
$\overline{\alpha}$ and the circular arc
$\{z=r_1e^{i\theta}:\,\theta_0\le \theta\le 2\pi-\theta_0\}$.  %
\end{rem} %
\smallskip

\begin{figure} 
\begin{minipage}{1.0\textwidth}
$$\includegraphics[scale=0.6,angle=0]{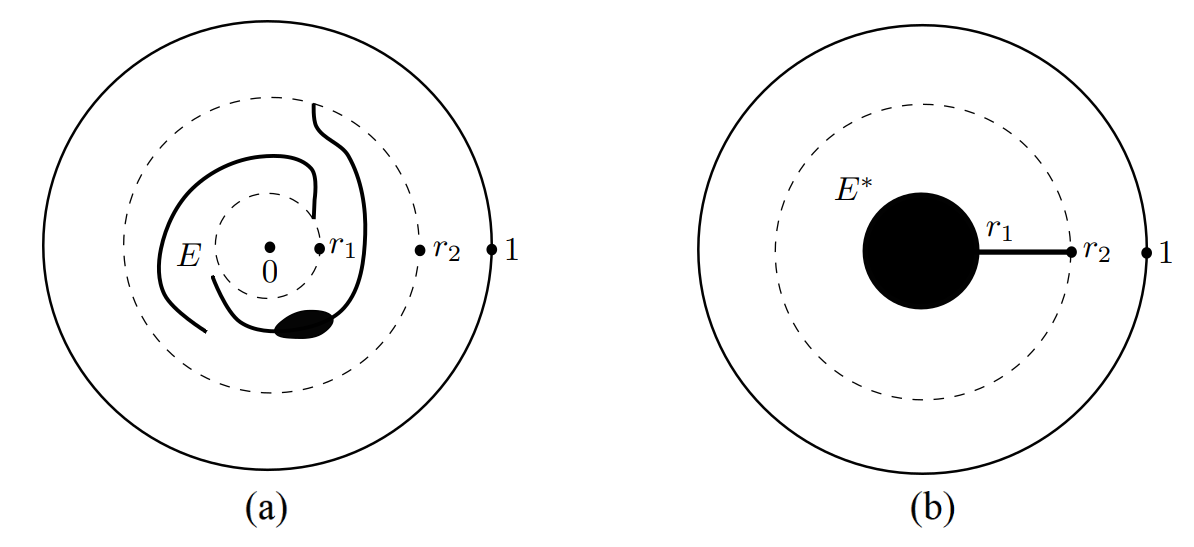}
$$
\end{minipage}
\caption{An admissible set and conjectured extremal set of Problem~\ref{Problem 5.25}.}%
\end{figure}


Now we turn to the hyperbolic version of the radial symmetrization
which in the Euclidean setting was introduced by G.~Szeg\"{o}
\cite{Szego}. To define this radial symmetrization, we need the
following notation. For $0<r<1$, $\alpha\in \mathbb{R}$, and a
compact set $E\subset \mathbb{D}$, define $E(r,\alpha)=E\cap
[re^{i\alpha},e^{i\alpha})$.

\begin{definition} \label{Definition 5.26}%
        Let $E$ be a compact set in $\mathbb{D}$ such that
        $\overline{\mathbb{D}}_r\subset E$, $0<r<1$. Then the radial
        symmetrization of $E$ with respect to $0$ is the compact set $E^{\rm rad}
        \subset \mathbb{D}$ with the property: for every $\alpha\in \mathbb{R}$, $E^{\rm rad}(r,\alpha)$ is a radial interval such that%
        \begin{equation} \label{Equation 5.24} %
        \int_{E^{\rm
rad}(r,\alpha)}\frac{|dz|}{|z|}=\int_{E(r,\alpha)}\frac{|dz|}{|z|}.
        \end{equation}
\end{definition} %

Equation(\ref{Equation 5.24}) shows that Szeg\"{o}'s radial
symmetrization is a rearrangement that is equimeasurable with
respect to the logarithmic metric. It found important applications
to several problems in Complex Analysis and Potential Theory.
However, because of the usage of the logarithmic metric, its
applications to the hedgehog problems studied in this paper are
rather limited. Indeed, certain compact sets on the interval
$[0,1]$  having a big hyperbolic length are transformed by this
symmetrization to the radial intervals with a very small
hyperbolic length.

In the search for  better estimates for the studied
characteristics of compact sets in $\mathbb{D}$, we turned to the
following version.

\begin{definition} \label{Definition 5.28}%
        The hyperbolic  radial symmetrization $E_h^{\rm rad}$ of a compact set
        $E\subset \mathbb{D}$ with respect to $0$ is defined to be a
        compact  set starlike with respect to $0$  and such that %
       $$ 
        \ell_{\mathbb{D}}(E_h^{\rm rad}\cap
        [0,e^{i\alpha}))=\ell_{\mathbb{D}}(E\cap
        [0,e^{i\alpha})) \quad {\mbox{for all $\alpha\in \mathbb{R}$.}} %
        $$ 
\end{definition}

        \medskip

        We recall that the compact sets in Theorems~\ref{Theorem 3.26}, \ref{Theorem 3.30}, and \ref{Theorem 3.42}
        having the minimal conformal capacity among all sets admissible
        for these theorems can be obtained via the hyperbolic radial
        symmetrization defined above. Also, the graphs in Figure~3 of the
        results of numerical computations
        suggest that the hyperbolic radial symmetrization of two radial
        intervals reduces the conformal capacity of these intervals.
        Therefore, the following conjecture sounds plausible.

        \begin{problem} \label{Problem 5.29}%
                Let $E$ be a compact set in $\mathbb{D}$. Prove (or disprove) that%
                $$ 
                {\rm cap}(E_h^{\rm rad})\le {\rm cap}(E).
                $$ 
        \end{problem} %

        We conclude this section with the following result, which
        describes how the hyperbolic area and hyperbolic diameter of a
        compact set $E$ behave under the hyperbolic radial symmetrization.

        \begin{lemma} \label{Lemma 5.30}%
                The hyperbolic area and the hyperbolic diameter do not increase under
                the hyperbolic radial symmetrization of $E$. %
        \end{lemma} %

        \noindent %
        \emph{Proof.}
        First we deal with the hyperbolic area $A_h(E)$ of a compact set $E$ in $\mathbb D$.
        Since every compact set in $\mathbb D$ can be approximated by a finite union
        of polar rectangles, we may assume that $E$ is a finite union of sets of the form
        $$
        R=\{re^{it}: r_1\leq r\leq r_2,\;t_1\leq t\leq t_2\},
        $$
        where $0\leq r_1<r_2<1$ and $0\leq t_1<t_2\leq 2\pi$.
        For such an $E$, the hyperbolic radial symmetrization  $E_h^{\rm rad}$ is again
        a union of sets of the same form, and moreover, each of the polar rectangles
        of $E_h^{\rm rad}$ is obtained by moving a rectangle of $E$ radially towards the origin,
        keeping its hyperbolic height $\uptau$ fixed. More precisely, if $R$ (as defined above)
        is a rectangle of $E$, then the corresponding rectangle of $E_h^{\rm rad}$ has the form
        $$
        \widehat{R}=\{re^{it}: \hat{r}_1\leq r\leq \hat{r}_2,\;t_1\leq t\leq t_2\},
        $$
        with $\hat{r}_1\leq r_1$ and
        \begin{equation}\label{Equation 5.31}
        \ell_{\mathbb D}([\hat{r}_1,\hat{r}_2])=\ell_{\mathbb D}([r_1,r_2])=\uptau.
        \end{equation}
        We express $r_2$ as a function of $r_1$ using (\ref{Equation 5.31}) and
        find that
        \begin{equation}\label{Equation 5.32}
        r_2=r_2(r_1)=(r_1+r(\uptau))/(1+r(\uptau)r_1),\quad {\mbox{where $r(\uptau)$
        is defined in (\ref{Equation 1.1}).}}
        \end{equation}
        Another elementary calculation gives
        $$ 
        A_h(R)=\frac{r_2^2-r_1^2}{2(1-r_1^2)(1-r_2^2)}(t_2-t_1).
        $$ 
        We consider the function
        $$
        g(r_1)=\frac{r_2^2-r_1^2}{(1-r_1^2)(1-r_2^2)}, \;\;r_2=r_2(r_1),\;\;0<r_1<1.
        $$
        By straightforward differentiation of $g(r_1)$ and taking into account (\ref{Equation 5.32}),
        we find that $g(r_1)$ is  strictly increasing. This implies that $A_h(R)\geq A_h(\widehat{R})$.
        Summing up over the hyperbolic areas of all rectangles, we conclude that
        $A_h(E)\geq A_h(E_h^r)$.

        \smallskip

        Next we turn to hyperbolic diameter, which we denote by ${\rm Diam}_h(E)$.
        We will use the following basic fact of the hyperbolic geometry.

        (a) Given a hyperbolic geodesic $\gamma$ and a point
        $a\not\in \gamma$, there is a unique hyperbolic geodesic
        $\gamma_\perp\ni a$ orthogonal to $\gamma$. Let
        $\gamma_\perp$ intersect $\gamma$ at $z=b$. Then the
        hyperbolic distance $d_\mathbb{D}(a,z)$ strictly increases
        as $z$ moves along $\gamma$ from $b$ to $\mathbb{T}$.

        To prove this result, we may assume that $\gamma=(-1,1)$
        and $a=is$, $0\le s<1$. Then $\gamma_\perp=(-i,i)$ and the
        monotonicity of $d_\mathbb{D}(is,r)$ for $0\le r<1$
        follows after simple calculations.

        Furthermore, if $\gamma=(-1,1)$ and $\Re a>0$, $\Im a>0$,
        then $\gamma_\perp$ intersects $\gamma$ at the point $b$
        such that $0<b<\Re a$.

        \smallskip

        We will also use the  following result that can be checked
        by standard Calculus technique.

        \smallskip

        (b) Let $z_1=re^{i\theta}$, $0<r<1$, and $z_2=e^{i\alpha}z_1$, $0\le
        \alpha\le \pi$. Then
        $p_\mathbb{D}(z_1,z_2)=\frac{2r\sin(\alpha/2)}{\sqrt{1-2r^2\cos\alpha+r^4}}$.
        Furthermore, the pseudo-hyperbolic distance
        $p_\mathbb{D}(z_1,z_2)$,
        and therefore the hyperbolic distance
        $d_\mathbb{D}(z_1,z_2)$, is a strictly increasing function
        of $r$, $0\le r<1$ and a strictly increasing function of
        $\alpha$, $0\le \alpha\le \pi$.

        \medskip

        Now let $z_1,z_2$ be two points on $E_h^{\rm rad}$ such that $d_{\mathbb D}(z_1,z_2)={\rm Diam}_h(E_h^{\rm rad})$.
        Note that $[0,z_1]\subset E_h^{\rm rad}$ and
        $[0,z_2]\subset E_h^{\rm rad}$. We may assume that $z_1=r_1\in (0,1)$, and $z_2=0$ or $z_2=r_2e^{it}$, $0<r_2\leq r_1$, $0\leq t\leq \pi$.
        If $z_2=-r_2$, $0\le r_2\le r_2$, then %
        $$
        {\rm Diam}_h(E_h^{\rm
        rad})=\ell_\mathbb{D}([-r_2,r_1])=\ell_\mathbb{D}(E\cap (-1,1))\le {\rm
        Diam}_h(E)%
        $$ %
        and the required result is proved.

        Next, we assume that $z_2=r_2e^{it}$ with $0<r_2\le r_1$,
        $0<t<\pi$. The set $E$ contains points $\zeta_1=s_1$ and $\zeta_2=s_2e^{it}$ with $s_1=\max\{|z|:\,z\in E\cap
        [0,1)\}$, $s_2=\max\{|z|:\,z\in E\cap [0,e^{it})\}$ such that
        $s_1\ge r_1$, $s_2\ge r_2$. If $\pi/2\le t <\pi$, then it
        follows from the monotonicity property stated in (a) that
        $d_\mathbb{D}(\zeta_1,\zeta_2)\ge d_\mathbb{D}(z_1,z_2)$
        and therefore ${\rm Diam}_h(E)\ge {\rm Diam}_h(E_h^{\rm
        rad})$ in this case.

        Now, we consider the case when $0<r_2\le r_1$,
        $0<t<\pi/2$. In this case, the hyperbolic geodesic
        $\gamma\ni z_1$ orthogonal to $(-e^{it},e^{it})$ crosses
        $(-e^{it},e^{it})$ at the point $z^*=r^*e^{it}$ with
        $0<r^*<r_1$. If $0<r_2\le r^*$, then
        $$ %
         {\rm Diam}_h(E_h^{\rm rad})=d_{\mathbb
D}(r_1,r_2e^{it})< d_{\mathbb D}(r_1,0)= \ell_{\mathbb D}(E_h^{\rm
rad}\cap [0,1))
        = \ell_{\mathbb D}(E\cap [0,1))\le  {\rm Diam}_h(E),
        $$
where the first inequality follows from the monotonicity property
stated in (a).

Thus, we are left with the case $r^*<r_2\le r_1$, $0<t<\pi/2$. We
work with the points $\zeta_1=s_1$ and $\zeta_2=s_2e^{it}$ defined
above. If $r^*<s_2\le r_1\le s_1$, then using  twice the
monotonicity property stated in (a) we obtain the required result:
 $$ %
         {\rm Diam}_h(E_h^{\rm rad})=d_{\mathbb
D}(z_1,z_2)< d_{\mathbb D}(z_1,\zeta_2)\le
d_\mathbb{D}(\zeta_1,\zeta_2)\le {\rm Diam}_h(E).
        $$

        If $r_1 <s_k\le  s_j$, $k\not=j$, then
 $$ %
         {\rm Diam}_h(E_h^{\rm rad})=d_{\mathbb
D}(z_1,z_2)\le d_{\mathbb D}(z_1,r_1e^{it})\le
d_\mathbb{D}(s_k,s_ke^{it})\le  d_\mathbb{D}(\zeta_k,\zeta_j)\le
{\rm Diam}_h(E),
        $$
where the first and the third inequalities follow from the
monotonicity property stated in (a) and the second inequality
follows from the monotonicity property with respect to $r$ stated
in (b).   Now, the inequality   ${\rm Diam}_h(E_h^{\rm rad})\le
{\rm Diam}_h(E)$ is proved in all cases. \hfill $\Box$

        \section{Hedgehog problems in ${\mathbb{R}}^3$}

        In this section we briefly mention how some of our results proved
        in the previous sections can be generalized for compact sets in
        the ball $\mathbb{B}=\{\overline{x}:\,|\overline{x}|<1\}$ in
        $\mathbb{R}^3$. Here, $\overline{x}=(x_1,x_2,x_3)$,
        $|\overline{x}|=\sqrt{x_1^2+x_2^2+x_3^2}$. The conformal capacity
        of a compact set $E\subset \mathbb{B}$ is defined as%
        $$ 
        {\rm{cap}}(E)=\inf \int_{\mathbb{B}} |\nabla u|^3\,dV,
        $$ 
        where $dV$ stands for the three-dimensional Lebesgue measure and
        the infimum is taken over all Lipschitz functions $u$ such
        that $u=0$ on $\partial  \mathbb{D}$ and $u=1$ on $E$. %

        The conformal capacity of $E$ is just the conformal capacity of
        the condenser $(\mathbb{B},E)$ and therefore it is invariant under
        M\"{o}bius automorphisms of $\mathbb{B}$, see \cite{G}. It is also
        known that the conformal capacity does not increase under polarization, see \cite{B}.
        However, it is an open problem to prove that the conformal
        capacity is strictly decreasing under polarization unless the
        result of the polarization and the original compact set coincide
        up to reflection with respect to the plane of polarization and up
        to a set of zero logarithmic capacity; see, for instance,
        \cite{Solynin2022}, where this question was discussed in the
        context of the Teichm\"{u}ller's problem in $\mathbb{R}^3$. This
        is why statements on the equality cases will be missing in all
        results, which we mention below.

        Furthermore, the contraction principle, as it is stated in
        Proposition~\ref{Proposition 2.19}, is not available in the context of the conformal
        capacity in dimensions $n\ge 3$. This issue was discussed, for
        instance, in Section~5.6.2 in \cite{B}.

        Below we list possible extensions of our results to
        $\mathbb{R}^3$ and to spaces of higher dimension.

        \begin{enumerate}
                \item[(1)] The non-strict monotonicity property in Lemma~\ref{Lemma 3.1} for
                a finite number of intervals on the diameter remains true in any
                dimension $n\ge 3$. %
                \item[(2)] The inequality (\ref{Equation 3.5}) of Lemma~\ref{Lemma 3.4} holds in any
                dimension $n\ge 3$.  %
                \item[(3)] The inequality (\ref{Equation 3.27}) of Theorem~\ref{Theorem 3.26} holds true for a
                set $E$ lying on $m$, $1\le m\le 6$, distinct radial intervals
                $I_k$ in $\mathbb{B}$ under the assumption that the angles between
                these intervals are greater than or equal to $\pi/2$.  %
                \item[(4)] Theorems~\ref{Theorem 3.30} and \ref{Theorem 3.42} also can be extended to the case
                of appropriate compact sets in $\mathbb{B}$.  %
                \item[(5)] The result of Theorem~\ref{Theorem 4.11} can be extended to the
                three-dimensional space. For this, under the assumptions of this
                theorem, we assume additionally that the unit disk $\mathbb{D}$ is
                embedded into $\mathbb{B}$ as follows: $\mathbb{D}=\mathbb{B}\cap
                \{(x_1,x_2,x_3)\in \mathbb{R}^3:\,x_3=0\}$. Then, inequality
                (\ref{Equation 4.12}) remain
                true for conformal capacities in $\mathbb{R}^3$.  %
        \end{enumerate}  %

        Some of the proofs in this paper make essential use of methods only
        available in the planar case and therefore the corresponding
        result for dimensions $n\ge 3$ remain open. Next we list a
        few of these open problems. %
        \begin{enumerate} %
                \item[(1)] Our proofs of Theorems~\ref{Theorem 3.35} and \ref{Theorem 4.1} require
                computations related to conformal mapping of doubly connected
                domains, which is not
                available in higher dimensions. %
                \item[(2)] The result stated in Theorem~\ref{Theorem 3.44} easily follows from
                the contractions principle. But, as we have mentioned above, the
                contraction principle is not available to problems on
                conformal capacity in dimensions $\ge 3$. %
                \item[(3)] Theorem~\ref{Theorem 4.8} can be proved by two methods: the first one
                uses the contraction principle and the second method uses explicit
                calculations. Both methods are not available in higher dimensions. %
                \item[(4)] To prove  Theorem~\ref{Theorem 4.14}  in the three-dimensional setting, we would need inequality
                (\ref{Equation 4.15}) under the assumption that the central body $E_0$
                is a ball $\overline{\mathbb{B}}_r=\{x\in \mathbb{R}^3:\,|x|\le r\}$, $0\le
                r<1$, and $E_1\subset [r,1)$ is the same as in
                Theorem~4.18.
                   The  proof of Theorem~\ref{Theorem 4.14} is based on
                 Baernstein's $*$-function method and the contraction principle.
                 Both techniques are still
                waiting to be developed for the case of this type of problems in $\mathbb{R}^3$.
                Thus, the counterpart of Theorem~\ref{Theorem 4.14} remains
                unproved in dimensions $n\ge 3$. %
        \end{enumerate} %

        \smallskip


\noindent%
 \textbf{Acknowledgements.} We dedicate this paper to the memory
of Jukka Sarvas whose work on symmetrization is a standard
reference in the potential-theoretic study of isoperimetric
inequalities and symmetrization.

We are indebted to Prof. M. M.S. Nasser and Dr. Harri Hakula for
kindly providing several numerical results for this paper. We also
thank K. Zarvalis for his help with the
figures.

We are grateful to the referee for many valuable suggestions.

\end{document}